\documentclass[12pt]{article}
\usepackage{graphicx} 
\usepackage{amsmath,amsthm}
\usepackage[a4paper,top=3cm,bottom=2cm,left=2cm,right=2cm,marginparwidth=1.75cm]{geometry}
\usepackage[backend=biber,sorting=none,style=numeric, citestyle=numeric-comp]{biblatex}
\newtheorem{definition}{Definition}%
\newtheorem{theorem}{Theorem}
\newtheorem{remark}{Remark}%
\newtheorem{lemma}{Lemma}%
\newtheorem{proposition}{Proposition}%
\usepackage[labelformat=simple]{subcaption}
\usepackage[export]{adjustbox}
\usepackage{lipsum}
\usepackage{bm}
\usepackage{bbm} 
\usepackage{amsfonts}
\usepackage{graphicx}
\usepackage{epstopdf}
\usepackage{algorithmic}
\usepackage{enumerate, enumitem}
\usepackage{comment}

\ifpdf
  \DeclareGraphicsExtensions{.eps,.pdf,.png,.jpg}
\else
  \DeclareGraphicsExtensions{.eps}
\fi

\title{Convex Approximations of Random Constrained Markov Decision Processes}
\usepackage{authblk}
\author[1]{V Varagapriya}
\author[1]{Vikas Vikram Singh}
\author[2]{Abdel Lisser}

\affil[1]{Department of Mathematics, Indian Institute of Technology Delhi, Hauz Khas, New Delhi, 110016, India\\
\texttt{varagapriyav@gmail.com}, \texttt{vikassingh@iitd.ac.in}}

\affil[2]{Laboratory of Signals and Systems, University Paris Saclay, CNRS, CentraleSupelec, Bat Breguet, 3 Rue Joliot Curie, Gif-sur-Yvette, 91190, France\\
\texttt{abdel.lisser@l2s.centralesupelec.fr}}
\date{}

\usepackage[dvipsnames, monochrome]{xcolor}
\usepackage{amsopn}

\usepackage{hyperref}
\usepackage{cleveref}
\crefname{equation}{}{} 
\crefname{theorem}{Theorem}{Theorem}

\crefrangelabelformat{subequation}{(#3#1#4-#5\crefstripprefix{#1}#2#6)} 
 \newcommand{\crefdefpart}[2]{%
  \hyperref[#2]{\namecref{#1}~\labelcref*{#1}.~\ref*{#2}}%
} 
\newtheorem{assumption}{Assumption}

\usepackage{booktabs, makecell, longtable}

\DeclareMathOperator*{\argmin}{arg\,min}

\iftrue
\newcommand{\comments}[1]{\footnote{\textcolor{blue}{\textit{#1}}}}
\newcommand{\AL}[1]{\textcolor{green}{#1}}
\newcommand{\VV}[1]{\textcolor{blue}{#1}}
\newcommand{\VS}[1]{\textcolor{red}{#1}}
\else
\newcommand{\comments}[1]{}
\newcommand{\AL}[1]{#1}
\newcommand{\VV}[1]{#1}
\newcommand{\VS}[1]{#1}
\fi
\begin{document}
 
\maketitle
\begin{abstract}
\noindent
Constrained Markov decision processes (CMDPs) are used as a decision-making framework to study the long-run performance of a stochastic system. It is well-known that a stationary optimal policy of a CMDP problem under discounted cost criterion can be obtained by solving a linear programming problem when running costs and transition probabilities are exactly known. In this paper, we consider a discounted cost CMDP problem where the running costs and transition probabilities are defined using random variables. Consequently, both the objective function and constraints become random. We use chance constraints to model these uncertainties and formulate the uncertain CMDP problem as a joint chance-constrained Markov decision process (JCCMDP).
Under random running costs, we assume that the dependency among random constraint vectors is driven by a Gumbel-Hougaard copula. 
Using standard probability inequalities, we construct convex upper bound approximations of the JCCMDP problem under certain conditions on random running costs.
In addition, we propose a linear programming problem whose optimal value gives a lower bound to the optimal value of the JCCMDP problem. 
 When both running costs and transition probabilities are random, we define the latter variables as a sum of their means and random perturbations. Under mild conditions on the random perturbations and random running costs, we construct convex upper and lower bound approximations of the JCCMDP problem.
We analyse the quality of the derived bounds through numerical experiments on a queueing control problem for random running costs.
 For the case when both running costs and transition probabilities are random, we choose randomly generated Markov decision problems called Garnets for numerical experiments. 
\end{abstract}

\noindent\textbf{Keywords:} 
Chance-constrained optimization, Constrained Markov decision process, Convex programming problem, Queueing control problem.

\section{Introduction}
A Markov decision process (MDP) is a mathematical framework used in the decision-making of a stochastically evolving system over time. It consists of a single decision-maker, a state space $S$, and an action set $A(s)$ available at each state $s\in S$. We assume $S$ and $A(s)$, $s\in S $, to be finite sets. 
Let the set of all state-action pairs be denoted by $\mathcal{K}$, i.e.,  $\mathcal{K}=\{(s, a) \mid s\in S, a \in A(s)\}$. We consider a system that evolves over discrete time points, where at time $t=0$, the system starts from a state $s_0$,   according to an initial distribution $\gamma$, i.e., the state $s_0$ is chosen with probability $\gamma(s_0)$. If the decision-maker chooses an action $a_0 \in A(s_0)$, a running cost $c(s_0, a_0)$ is incurred, and at the next time point, the system moves to another state $s_1$ with a transition probability $p(s_1 \vert s_0, a_0)$.  This process repeats infinitely. 
The aim of the decision-maker is to minimize the overall cost incurred. In many real-world instances, however, the decision-maker incurs multiple running costs each time  \cite{altman,puterman}.
For instance, if at $t=0$, the decision-maker chooses an action $a_0 \in A(s_0)$,  running costs $c(s_0, a_0)$ and $d^k(s_0, a_0)$,  $k \in \mathbb{K}$, are incurred; $\mathbb{K} = \{ 1,2,\ldots, K \}$. In such cases, the aim is to minimize the overall cost associated with $c$  such that
each overall cost associated with $d^k$, $k\in \mathbb{K}$, is upper bounded by a constant $\xi_k$. Such a class of MDP is called a constrained Markov decision process (CMDP).

At each time point $t$, the decision-maker may choose an action according to a decision rule which depends on the history $h_t$ of state-action pairs traversed till time $t-1$  and a state $s_t$ at time $t$, i.e., $h_t = (s_0, a_0, s_1, a_1,\ldots,s_{t-1}, a_{t-1},s_t)$. Such a decision rule is called history-dependent and it is denoted by $f_t^h=f_t(h_t)\in \wp(A(s_t))$, where $\wp(A(s_t))$ denotes the set of probability distributions over the action set $A(s_t)$. In contrast, when the choice of each action depends only on the state (i.e., it is independent of history and time), the decision rule is called stationary and denoted by $f$. A sequence of history-dependent decision rules  $f^h=(f_t^h)_{t=0}^\infty$ is called a history-dependent policy, while a sequence of stationary decision rules  $f$ (denoted with abuse in notation) is called a stationary policy. A stationary policy $f$ can be viewed as a vector $\big(f(s)\big)_{s\in S}$, where $f(s)\in \wp(A(s))$ for each $s\in S$.  We denote the set of all history-dependent and stationary policies by $F_{HD}$ and $F_S$, respectively. 
It is well-known that there exists a stationary deterministic optimal policy for an MDP problem and it can be computed using variants of dynamic programming-based methods as well as linear programming (LP)-based methods \cite{puterman,bertsekas2005dynamic}. On the other hand, for a CMDP problem, there exists a stationary randomized optimal policy which can be computed using an LP problem \cite{altman}.

Let $(\mathbb{X}_t,\mathbb{A}_t)$ denote the state-action pair at time $t$. For a given $f^h \in F_{HD}$ and an initial distribution $\gamma$, let $\mathbb{P}_{\gamma}^{f^h}$ denote the probability measure over the state-action trajectories  while  $\mathbb{E}_{\gamma}^{f^h}$ denote the expectation operator corresponding to $\mathbb{P}_{\gamma}^{f^h}$ (see Section 2.1.6 of \cite{puterman}). In this paper,  the overall costs are considered to be  expected discounted  costs with a discount factor $\alpha \in (0,1)$, and are defined as 
\begin{align*}
     C_{\alpha}(\gamma,f^h) & = (1-\alpha) 
\sum_{t=0}^{\infty}\alpha^{t}\mathbb{E}_{\gamma}^{f^h}c(\mathbb{X}_t,\mathbb{A}_t),
\nonumber \\
   D^k_\alpha(\gamma,f^h) & = (1-\alpha)\sum_{t=0}^\infty \alpha^t  \mathbb{E}_{\gamma}^{f^h}d^k(\mathbb{X}_t,\mathbb{A}_t), \ \forall \ k \in \mathbb{K}.
  \end{align*}
An optimal policy of a CMDP problem is obtained by solving the following optimization problem
\begin{align}\label{COP_HD_policy}
 \min_{f^h \in F_{HD}} & \  \ C_\alpha(\gamma,f^h) \nonumber \\
  \text{s.t.} & \  \ D^k_\alpha(\gamma,f^h) \leq \xi_k, \ \forall \ k \in \mathbb{K}.
\end{align}
From Theorem 3.1 of \cite{altman}, \eqref{COP_HD_policy} can be restricted to the class of stationary policies without loss of optimality and from Theorem 3.2 of \cite{altman}, it can be equivalently written as the following LP problem.
\begin{align}\label{equi_LP}
    \min_{\rho\in \mathbb{R}^{|\mathcal{K}|}} & \  \sum_{(s,a) \in \mathcal{K}} \rho(s,a)c(s,a) \nonumber\\
    \text{s.t.} & \  \sum_{(s,a) \in \mathcal{K}} \rho(s,a)d^k(s,a)\leq \xi_k, \ \forall \ k \in \mathbb{K}, \  \rho \in \mathcal{Q}^{\alpha}(\gamma), 
\end{align}
where
\begin{align}\label{definition_q_alpha}
\mathcal{Q}^{\alpha}(\gamma) = \Big\{\rho\in \mathbb{R}^{|\mathcal{K}|} \ \bigl \lvert \ \sum_{(s,a) \in \mathcal{K}} \rho(s,a) \big(\delta({s}', s)
 -\alpha p({s}'|s,a)\big) 
=(1-\alpha) \gamma({s}'),\ \forall \ {s}' \in S,  \nonumber\\  \rho(s,a)\geq 0, \ \forall \ (s,a) \in \mathcal{K}\Big\},
\end{align}
and $\delta({s}', s)$ denotes the Kronecker delta. As a result, if $\big( \rho^*(s,a) \big)_{(s,a) \in \mathcal{K}}$ is an optimal solution of \eqref{equi_LP}, then an optimal stationary policy $f^*$ of \eqref{COP_HD_policy} can be obtained by the  relation $ \displaystyle 
f^*(s,a) = \displaystyle \frac{\rho^*(s,a)}{\displaystyle\sum_{a \in A(s)}\rho^*(s,a)}$,    $(s,a) \in \mathcal{K}  
 $, if the denominator is non-zero. Otherwise, we select an arbitrary $f(s)\in \wp(A(s))$.

 The reformulation of \eqref{COP_HD_policy} into \eqref{equi_LP} is derived assuming that the model parameters, namely, the running costs and the transition probabilities, do not vary with time, i.e., they are stationary and exactly known. The former assumption is made in a majority of prior work (Chapter 6 of \cite{puterman}, Chapter 7 of \cite{bertsekas2005dynamic}), although the theory of MDP problems is also developed for non-stationary model parameters. For instance,  in \cite{ghate2013linear,cheevaprawatdomrong2007solution},  the authors considered non-homogeneous MDP problems and discussed specialized algorithms to compute an optimal policy. In contrast, the latter assumption poses a considerable risk of resulting in erroneous policies, since the model parameters are typically estimated using historical data or learnt with experience \cite{mannor2007bias}. Several studies have discussed MDP problems with uncertain model parameters using a robust optimization framework
\cite{satia1973, givan2000bounded, nilim2005robust, iyengar2005robust, wiesemann2013robust}. In \cite{cmdp_cost_uncertainty,Rankone_TPM}, a robust CMDP problem under uncertain running costs and transition probabilities are considered.    However, in \cite{erick}, the authors argued that this framework leads to a conservative solution and used a chance-constrained optimization as an alternative framework. They considered MDP problems where either running costs or transition probabilities are random and constructed a chance-constrained MDP problem in each case. 
They showed that this problem is equivalent to a second-order cone programming (SOCP) problem
when the running costs follow a Gaussian distribution and the transition probabilities are exactly known. On the other hand, for the case of known costs and Dirichlet distributed transition probabilities, approximations are proposed. Varagapriya et al. \cite{varagapriya2022joint} extended this theory to a CMDP problem by formulating it as a joint chance-constrained Markov decision process  (JCCMDP) when the running costs are random and transition probabilities are known.  When the running cost vectors follow elliptically symmetric distribution and the dependency among the random constraints is driven by a Gumbel-Hougaard copula, the authors proposed two SOCP-based approximations whose optimal values provide upper and lower bounds to the optimal value of the JCCMDP problem. In general, such a model is severely limited in practice, since it is assumed that the exact distribution of the model parameters is known. As a result, unlike in \cite{erick, varagapriya2022joint},  we make no such assumptions in this paper. We consider a JCCMDP problem where either running costs or both running costs and transition probabilities are random. To the best of our knowledge, this is the first paper that considers the case when the running costs and transition probabilities together are random in MDPs/CMDPs. 
  Due to the presence of a joint chance constraint, the JCCMDP problem is, in general, non-convex and exceedingly complex to solve. A natural 
approach is to construct convex approximations which provide upper and lower bounds to the optimal value of the JCCMDP problem.    From the optimal solution of the upper bound approximation, we obtain an $\epsilon$-optimal policy, where $\epsilon$ represents the difference between the upper and lower bounds. For a sufficiently small $\epsilon$ value, a $\epsilon$-optimal policy can be considered as good as an optimal policy. 
 The JCCMDP problem under random running costs generalizes the results from \cite{varagapriya2022joint}. Additionally, we study the JCCMDP problem when running costs and transition probabilities together are random, and propose novel convex approximations.
 The contributions of the paper are summarized as follows.
 \begin{enumerate}[wide]
     \item  
     We consider a JCCMDP problem under random running costs and known transition probabilities, where the dependency among the random constraint vectors is driven by a Gumbel-Hougaard copula. Unlike in \cite{varagapriya2022joint},
     the running cost vectors need not have a known distribution.
     Using standard probability inequalities, we construct four different convex upper bound  approximations of the JCCMDP problem. Each approximation assumes a different set of conditions on the random running cost vectors. 
   Furthermore, we show that a lower bound to the optimal value of the JCCMDP problem can be obtained by solving an LP problem.
     \item We consider a JCCMDP problem under both random running costs and transition probabilities. The transition probabilities are defined as a sum of their means and random perturbations. Similar to the previous case, we construct three convex upper bound approximations using standard probability inequalities. Each approximation assumes a different set of conditions on the random perturbations and random running cost vectors. 
    Furthermore, we show that a lower bound to the optimal value of the JCCMDP problem can be obtained by solving an LP problem.  
     \item  To analyse the gap between the upper and lower bounds obtained from the proposed approximations, we perform numerical experiments on randomly generated instances of a queueing control problem, as well as on a class of MDP problems called  Garnets. 
 \end{enumerate}
 The structure of the rest of the paper is as follows. Section \ref{model} contains the definition of the JCCMDP problem and preliminary results on standard probability inequalities. In Sections \ref{JCCMDP problem under random running costs} and  \ref{JCCMDP problem under random running costs and transition probabilities}, we present convex upper and lower bound approximations when only running costs are random and when both running costs and transition probabilities are random, respectively. Section \ref{Num} contains numerical experiments on randomly generated instances. We conclude the paper in Section~\ref{sec:conclusions}.

 \section{Joint chance-constrained Markov decision processes} \label{model}
 In this paper,  CMDP problems under random running costs and random transition probabilities are examined. We take into account the case where either running costs or both running costs and transition probabilities are random and examine them separately. The following is a summary of the notations used in each case.
 \begin{enumerate}[wide,label=C.{\arabic*}.]
    \item \label{JCCMDP_running_costs_desc.} We assume that the transition probabilities are known and stationary, and are denoted by $p(s'|s, a)$,   $(s, a)\in \mathcal{K}$, $s'\in S$.  The running costs are random and stationary, and are defined on a probability space $(\Omega,\mathcal{F},\mathbb{P})$. For a given realization  $\omega \in \Omega$, let $\tilde{c}(s,a,\omega)$ and $\tilde{d}^k(s,a,\omega)$,  $k \in \mathbb{K}$, denote the running costs incurred at the state-action pair $(s,a) \in \mathcal{K}$. For a given $f^h \in F_{HD}$ and initial distribution $\gamma$, we denote the corresponding random expected discounted costs by $\tilde{C}_{rc}(\gamma,f^h)$ and $\tilde{D}^k_{rc}(\gamma,f^h)$,   $k \in \mathbb{K}$. 
    \item \label{JCCMDP_transition_probabilities_desc.} We assume running costs and transition probabilities are random and stationary, and are defined on a probability space $(\Omega,\mathcal{F},\mathbb{P})$. For a given realization  $\omega \in \Omega$,  let $\tilde{p}(s' \vert s,a)(\omega)$ denote the transition probability from state-action pair $(s,a) \in \mathcal{K}$ to a state $s' \in S$. For a given $f^h \in F_{HD}$ and initial distribution $\gamma$, we denote the corresponding random expected discounted costs  by $\tilde{C}_{rc{\text-}tp}(\gamma,f^h)$ and $\tilde{D}^k_{rc{\text-}tp}(\gamma,f^h)$,  $k \in \mathbb{K}$.
\end{enumerate}
We consider the case where the aim of the decision-maker is to minimize the expected discounted cost in the objective function,  incurred with at least a given probability level $p_0$, such that the expected discounted costs in the constraints are jointly satisfied with at least a given   probability level $p_1$. This leads to the following JCCMDP problem
\begin{align}\label{common_JCCMDP}
         & \min_{z,f^h \in F_{HD}} \,  z \nonumber \\
        \text{s.t. }& \mathbb{P}\Big( \tilde{C}_\mathfrak{v}(\gamma,f^h)\leq z \Big)\geq p_0, \nonumber \\
         & \mathbb{P}\Big( \tilde{D}^k_\mathfrak{v}(\gamma,f^h)\leq \xi_k,  \ \forall \ k \in \mathbb{K} \Big)\geq p_1,
\end{align}
where $\mathfrak{v} \in \{rc, rc{\text-}tp \}$. When $\mathfrak{v} = rc$, \eqref{common_JCCMDP} becomes a JCCMDP problem under random running costs and represents the case C.1. On the other hand, when $\mathfrak{v} = rc{\text-}tp$, \eqref{common_JCCMDP} becomes a JCCMDP problem under both
random running costs and transition probabilities and represents the case C.2.
 We use inner approximations of a chance constraint based on standard probability inequalities to construct convex upper bound approximations of \eqref{common_JCCMDP}. 
 Thus,  the optimal value of each approximation provides an upper bound to the optimal value of \eqref{common_JCCMDP}.  In addition, we propose a new linear lower bound approximation of  \eqref{common_JCCMDP} such that its optimal value provides a lower bound to the optimal value of \eqref{common_JCCMDP}. In the following sections,  we present basic definitions and existing results on convex inner approximations of a general linear chance constraint. Furthermore, we propose its linear outer approximation. These approximations are used later in Sections \ref{JCCMDP problem under random running costs} and  \ref{JCCMDP problem under random running costs and transition probabilities} to derive subsequent results for the JCCMDP problems.

 \subsection{Basic definitions and results}\label{Preliminaries_with_common_notation}
 \begin{definition}[Copula \cite{Nelsen,Jaworski}]
     A copula of dimension $K$, where $K\geq 2$, is a distribution function defined on $[0,1]^K$,   whose one-dimensional marginals are uniformly distributed on $[0,1]$.
\end{definition}
In this paper, we consider a specific class of copula, namely, the Gumbel-Hougaard copula, denoted by $\mathcal{C}_{\theta}$, defined below.
\begin{definition}[Gumbel-Hougaard copula \cite{Nelsen,Jaworski, cheng2015chance}] \label{definition_GH_copula}
Let $u=(u_1,u_2,\dots, \allowbreak u_K)^T  \in [0,1]^K$ and $\exp$ denotes the exponential function. Then 
\begin{align}
     \mathcal{C}_\theta(u)=\exp\Bigg\{-\bigg[\sum_{k=1}^K \big(-\ln(u_k)\big)^{\theta}  \bigg]^{\frac{1}{\theta}}\Bigg\},
\end{align}
   where $\theta \geq 1$. 
\end{definition}
When $\theta = 1$, $\mathcal{C}_1(u) = \prod_{k = 1}^K u_k$ denotes the product copula and it corresponds to the case when $X_1, X_2,\ldots, X_K$ are independent.
For a given random vector $X = (X_1, X_2, \ldots, X_K)$, 
the relation between its joint distribution and a copula is given by Sklar's theorem stated below.
\begin{theorem}[Sklar's theorem \cite{Nelsen}] \label{Sklar_theorem}
Let $\hat{\Phi}$ be a given $K$-dimensional distribution function and $\hat{\Phi}_{1}, \hat{\Phi}_{2},\dots, \hat{\Phi}_{K}$ be its one dimensional marginals. Then there exists a copula $\mathcal{C}$ satisfying
\begin{align*}
    \hat{\Phi}(v_1,v_2,\ldots,v_K)=\mathcal{C}\big(\hat{\Phi}_1(v_1),\hat{\Phi}_2(v_2),\ldots,\hat{\Phi}_K(v_K)\big), \,\forall \,v_1,v_2,\ldots,v_K \in \mathbb{R}.
\end{align*}
If all the marginals are continuous, then $\mathcal{C}$ is unique throughout its domain; otherwise it is unique on 
$ \times_{i=1}^K$Range\,$\hat{\Phi}_i$,
where $Range\,\hat{\Phi}_i$, $i=1,2,\ldots,K$, denotes the set of all possible values $\hat{\Phi}_i$ can take.
\end{theorem}

\subsection{Approximations of a linear chance constraint}
Let $\mathfrak{z} = ( \mathfrak{z}_i )_{i = 1}^n $ be a random vector defined on a probability space $(\bf{\Omega}, \mathcal{F}, \mathbb{P})$,  $r = (r_i)_{i = 1}^n$ be an arbitrarily fixed vector, and $a$ be some fixed constant. We consider the following chance constraint
\begin{align}\label{general_chance_constraint} 
     \mathbb{P}\big( r^T \mathfrak{z}  \leq a  \big) \geq p,     
 \end{align}
 where $p \in (0,1)$. In  Theorems  \ref{C_H_Sub_G_common_statement_general_CC} and   \ref{Bernstein inequality}, we present existing 
 convex inner approximations of \eqref{general_chance_constraint}, including those  based on standard probability inequalities.
 \begin{theorem}\label{C_H_Sub_G_common_statement_general_CC}
 Let $f: (0,1) \rightarrow \mathbb{R}$ and $V:\mathbb{R} \rightarrow \mathbb{R}$. Consider the following general constraint 
 \begin{align}\label{common_C_H_sub_G_lin_CC}
    r^T \mu_{ \mathfrak{z} }  + f(p)  V(r) \leq a.
\end{align}
 \begin{enumerate}[wide]
\item  \label{One-sided Chebyshev inequality}
\textnormal{(One-sided Chebyshev inequality \cite{pinter1989deterministic})}
Let the mean vector and the covariance matrix of $\mathfrak{z}$, denoted by 
$\mu_{\mathfrak{z}} \in \mathbb{R}^{n  }$ and $\Sigma_{\mathfrak{z}} \in \mathbb{R}^{n \times n }$, respectively, be known.   Then, a feasible solution $r$ of the second-order cone constraint \eqref{common_C_H_sub_G_lin_CC}, when  $ \displaystyle f(p) = \sqrt{\frac{p} {1-p}}$ and $V(r) = \| \Sigma_{\mathfrak{z}}^{\frac{1}{2}} r \|_2$, is a feasible solution of  \eqref{general_chance_constraint}.  
\item \label{Hoeffding inequality}
\textnormal{(Hoeffding inequality \cite{Hoeffding1963inequalities})}
Let $\mu_{\mathfrak{z}}$, and  upper and lower bound vectors of $\mathfrak{z}$, denoted by     $\mathfrak{z}^u = (\mathfrak{z}_i^u)_{i = 1}^n$ and $\mathfrak{z}^l = (\mathfrak{z}_i^l)_{i = 1}^n$, respectively, be known, and the components of $\mathfrak{z} $ be independent. Then, a feasible solution $r \VV{\geq 0}$ of the second-order cone constraint \eqref{common_C_H_sub_G_lin_CC}, when $\displaystyle f(p) =  \sqrt{-\frac{1}{2} \ln (1-p)}$ and $V(r) = \| D_{ \mathfrak{z} }^{\frac{1}{2} } r \|_2$; $D_{\mathfrak{z}} $ is an $n \times n$ diagonal matrix whose $i^{th}$ diagonal entry is defined as $  (\mathfrak{z}_i^u - \mathfrak{z}_i^l)^2$,   is a feasible solution of \eqref{general_chance_constraint}.
\item \label{Sub-Gaussian condition}
\textnormal{(\!\cite{wainwright2019high})}
Let $\mu_{\mathfrak{z}} $ be known, and
   the components of $\mathfrak{z} $ be independent  sub-Gaussian random variables  such that $\mathfrak{z_i}$ has  sub-Gaussian parameter $\sigma_{\mathfrak{z}_i}$. Then, a feasible solution $r$ of the second-order cone constraint \eqref{common_C_H_sub_G_lin_CC}, when $f(p) = \sqrt{-2 \ln (1-p)}$ and $V(r) =\displaystyle \sqrt{ 
\sum_{i=1}^n r_i^2 \sigma_{\mathfrak{z}_i}^2   }$, is a feasible solution of \eqref{general_chance_constraint}. 
\end{enumerate}
\end{theorem}
\begin{remark}
    When  the conditions of part \ref{Hoeffding inequality} of  \Cref{C_H_Sub_G_common_statement_general_CC} holds,  $\mathfrak{z_i}$ is sub-Gaussian with parameter at most $\displaystyle \frac{\mathfrak{z}_i^u - \mathfrak{z}_i^l}{2}$, for each $i$ \cite{wainwright2019high}.
\end{remark}
\begin{theorem}[Bernstein inequality \cite{pinter1989deterministic}]\label{Bernstein inequality}
    Let $\mu_{\mathfrak{z}}$, $\mathfrak{z}^u$, and $\mathfrak{z}^l$  be known and the components of $\mathfrak{z} $ be independent. Then, a feasible solution $r \VV{\geq 0}$ of the  following convex  constraint 
    \begin{align*} 
        \sum_{i = 1}^n
        \ln \bigg\{  \bigg( \frac{ \mu_{\mathfrak{z}_i } - \mathfrak{z}_i^l}{\mathfrak{z}_i^u - \mathfrak{z}_i^l } \bigg)  \exp \big( h r_i \mathfrak{z}_i^u \big) + 
        \bigg( \frac{ \mathfrak{z}_i^u - \mu_{\mathfrak{z}_i } }{\mathfrak{z}_i^u - \mathfrak{z}_i^l } \bigg) \exp \big(h r_i \mathfrak{z}_i^l  \big) \bigg\} 
        \leq 
        \ln (1-p) + h a , 
    \end{align*}
     where $h > 0$ is an arbitrarily fixed constant,   is a feasible solution of  \eqref{general_chance_constraint}.
\end{theorem}
We propose a new outer approximation of \eqref{general_chance_constraint}, consisting of a set of linear inequalities. 
\begin{theorem}\label{general lower bound}
    Let  $\mu_{\mathfrak{z}}$, $\mathfrak{z}^u $, and  $\mathfrak{z}^l $ be known and $r \geq 0$ be a feasible solution of  \eqref{general_chance_constraint}. Then there exists an $\mathfrak{m}$ such that $(r, \mathfrak{m})$  is a feasible solution of the following set of linear constraints
    \begin{subequations}\label{general_lowerbound}
    \begin{align}
       &     r^T \mu_{ \mathfrak{z} } - a \leq (1-p) \mathfrak{m} , \  \lambda \leq \mathfrak{m}, \  r^T \mathfrak{z}^u - a \leq \mathfrak{m},  \label{general_lowerbound_expectation}  \\ 
       &  r^T \mathfrak{z}^l \leq a,\label{general_lowerbound_using_quantile}
       \end{align}
       where $\lambda > 0 $ is a constant. 
    \end{subequations}
\end{theorem}
\begin{proof}
    The chance constraint \eqref{general_chance_constraint} is the same as 
\begin{align}\label{general_chance_constraint_with_expectation}
        1-p \geq \mathbb{E}\big( \mathbbm{1}_{  \{ r^T \mathfrak{z} - a > 0  \} }\big),  
    \end{align}
where $\mathbbm{1}_{\{\cdot\}}$ and $\mathbb{E}(\cdot)$ denote the indicator function and the corresponding expectation, respectively. Since $r \geq 0$, 
for each realization $\omega \in \Omega$, we have
\begin{align*}
   &  \mathbbm{1}_{  \{ r^T \mathfrak{z}(\omega)  - a > 0  \}} \geq   0, \\ 
    &  \mathbbm{1}_{  \{ r^T \mathfrak{z}(\omega)  - a > 0  \}} \geq  \frac{r^T \mathfrak{z}(\omega) - a}{\mathfrak{m}}, 
\end{align*}
where   $ \mathfrak{m} = \max\big( \lambda,   r^T {\mathfrak{z}^u}  - a \big)$, for  some small $\lambda > 0$. Therefore, we have 
\begin{align*}
        \mathbb{E}\big( \text{$\mathbbm{1}$}_{  \{ r^T \mathfrak{z} - a > 0  \} }\big) \geq \max\bigg( 0, \frac{ r^T \mu_{ \mathfrak{z} } - a }{ \mathfrak{m} } \bigg). 
       \end{align*}
     Using \eqref{general_chance_constraint_with_expectation}, we obtain  
       \[
       \max\bigg( 0, \frac{ r^T \mu_{ \mathfrak{z} } - a }{ \mathfrak{m} } \bigg)\le 1-p.
       \]
    This implies that $(r, \mathfrak{m})$ is a feasible solution of \eqref{general_lowerbound_expectation}. 
    Furthermore, \eqref{general_chance_constraint} implies   $\Phi_{r^T \mathfrak{z}} (a) \geq p$ \VV{and} 
   $r\ge 0$ \VV{implies} $r^T \mathfrak{z}\ge   r^T\mathfrak{z}^l$. Therefore, $r$ satisfies \eqref{general_lowerbound_using_quantile}\VV{, since otherwise, $\Phi_{r^T \mathfrak{z}} (a) = 0$. Thus,}  the result follows.
\end{proof}

\section{JCCMDP problem under random running costs}\label{JCCMDP problem under random running costs}
In this section, we consider the JCCMDP problem for the case \ref{JCCMDP_running_costs_desc.}, i.e., the transition probabilities are known while the running costs are random. Varagapriya et al. \cite{varagapriya2022joint}  showed that \eqref{common_JCCMDP}, when $\mathfrak{v} = rc$, can be equivalently reformulated into an LP problem with a joint chance constraint. We present this result in the following lemma.
\begin{lemma}[Proposition 1 of \cite{varagapriya2022joint}] \label{lemma}
 The   JCCMDP problem
    \eqref{common_JCCMDP} when $\mathfrak{v} = rc$ can be restricted to the class of stationary policies without loss of optimality. Furthermore, we can equivalently write it as the following LP problem with a joint chance constraint   \begin{subequations}\label{equivalent_JCCMDP_random_costs_using_rho}
   \begin{align}
       \min_{z,\, \rho } & \ z \nonumber \\
        \textnormal{s.t.} 
 & \ 
\mathbb{P} \Bigg( \sum_{(s,a) \in \mathcal{K}}  \rho(s,a) \tilde{c}(s,a) \leq z \Bigg) \geq p_0, \label{equivalent_JCCMDP_random_costs_using_rho_individual_constraint}  \\
& \ 
\mathbb{P} \Bigg( \sum_{(s,a) \in \mathcal{K}}  \rho(s,a) \tilde{d}^k(s,a) \leq \xi_k, \ \forall \  k \in \mathbb{K} \Bigg) \geq p_1, \label{equivalent_JCCMDP_random_costs_using_rho_joint_constraint}  \\
& \ 
\rho \in \mathcal{Q}^{\alpha}(\gamma). 
   \end{align}
   \end{subequations}
     An optimal stationary policy,  $f^*$ of \eqref{common_JCCMDP} can be obtained from  
 an optimal solution $(z^*, \rho^*)$ of \eqref{equivalent_JCCMDP_random_costs_using_rho} by  $ \displaystyle f^*(s,a) = \frac{\rho^*(s,a)}{\sum_{a\in A(s)} \rho^*(s,a)}$,  $(s,a) \in \mathcal{K}$, if the denominator is non-zero. We arbitrarily choose  $f(s)\in \wp(A(s))$ otherwise.
\end{lemma}
The constraints \crefrange{equivalent_JCCMDP_random_costs_using_rho_individual_constraint}{equivalent_JCCMDP_random_costs_using_rho_joint_constraint} pose major difficulties in solving \eqref{equivalent_JCCMDP_random_costs_using_rho}. For a given $\rho \in \mathcal{Q}^{\alpha}(\gamma)$, checking the feasibility of  \Crefrange{equivalent_JCCMDP_random_costs_using_rho_individual_constraint}{equivalent_JCCMDP_random_costs_using_rho_joint_constraint} is NP-hard because it requires solving multi-dimensional integration. In addition, the feasible region need not be convex  \cite{luedtke2008sample, ahmed2014convex}. In order to alleviate these problems, 
we assume that the dependence among the constraint vectors is driven 
by a Gumbel-Hougaard copula whose separable nature enables us to equivalently transform the
joint chance constraint   \eqref{equivalent_JCCMDP_random_costs_using_rho_joint_constraint}
 into a set of individual chance constraints. In particular, we make the following assumption.
\begin{assumption}\label{dependence_joint_chance_constr}
    The dependence among the random variables in \eqref{equivalent_JCCMDP_random_costs_using_rho_joint_constraint} is driven by a Gumbel-Hougaard copula, given in \eqref{definition_GH_copula}, which is independent of the decision vector $\rho$.
 \end{assumption} 
 Under   \Cref{dependence_joint_chance_constr}, \eqref{equivalent_JCCMDP_random_costs_using_rho} can be equivalently written as an individual chance-constrained programming problem. We state it as a lemma whose proof is given in \cite{varagapriya2022joint}.  
 \begin{lemma}
    Under \Cref{dependence_joint_chance_constr}, the   problem \eqref{equivalent_JCCMDP_random_costs_using_rho} can be equivalently written as the following individual chance-constrained programming problem 
    \begin{subequations}\label{ICCP_random_costs}
    \begin{align}
        \min_{z, \rho, (y_k)_{k \in \mathbb{K} } } & \ z \nonumber \\
        \textnormal{s.t.} 
 & \ 
 \mathbb{P} \left(   \rho^T \tilde{c}  \leq z \right)
\geq p_0, \label{ICCP_random_costs_objective}   \\
& \ 
\mathbb{P} \left(   \rho^T \tilde{d}^k  \leq \xi_k  \right)
\geq p_1^{ {y}_k^{\frac{1}{\theta}} },  \ \forall \ k \in \mathbb{K},  \label{ICCP_random_costs_constraints}  \\
& \
\rho \in \mathcal{Q}^{\alpha}(\gamma),  \ \sum_{k \in \mathbb{K} } y_k = 1, \  y_k \geq 0,  \ \forall \ k \in \mathbb{K}. \label{ICCP_random_costs_sum_1_non_negative_flow_balance} 
    \end{align}
    \end{subequations}
\end{lemma}
\begin{proof}
    The proof follows from Proposition 2 of \cite{varagapriya2022joint}.
\end{proof}
Unlike in \cite{varagapriya2022joint}, we consider the case where the running cost vectors $\tilde{c}$ and $\tilde{d}^k$ need not follow elliptically symmetric distributions, rather, their distributions are unknown. 
Therefore, we use  Theorems \ref{C_H_Sub_G_common_statement_general_CC},  \ref{Bernstein inequality}, and \ref{general lower bound} to construct convex approximations whose optimal values provide upper and lower bounds to the optimal value of the JCCMDP problem.  We propose four upper bound approximations and one lower bound approximation. Each approximation assumes a subset of conditions given in Assumption \ref{assumption_costs_bounds_mean_known_independence} below. 
\begin{assumption}\label{assumption_costs_bounds_mean_known_independence} 
\begin{enumerate}[wide=0pt]
\item \label{Mean_Cov_RC_known} 
  The mean  vectors of  $\tilde{c}$ and $\tilde{d}^k$, denoted by  $\mu_{\tilde{c}} {=} \big( \mu_{\tilde{c}}(s, a) \big)_{(s, a) \in \mathcal{K}}$ and   $\mu_{\tilde{d}^k}= \big( \mu_{\tilde{d}^k}(s, \allowbreak a) \big)_{(s, a) \in \mathcal{K}}$, and   covariance matrices of $\tilde{c}$ and $\tilde{d}^k$, denoted by  $\Sigma_{\tilde{c}} \in \mathbb{R}^{\vert \mathcal{K} \vert \times \vert \mathcal{K} \vert} $ and   $\Sigma_{\tilde{d}^k} \in \mathbb{R}^{\vert \mathcal{K} \vert \times \vert \mathcal{K} \vert}$,   respectively, are known. 
  \end{enumerate}
\begin{enumerate}[wide,resume]
\item   \label{mean_UB_
LB_RC_Obj_Known} 
The vector $\mu_{\tilde{c}}$, and  upper and lower bound vectors of   $\tilde{c}$, denoted by 
$\tilde{c}^u = \big( \tilde{c}^u(s, a) \big)_{(s, a) \in \mathcal{K}}$ and 
$\tilde{c}^l= \big( \tilde{c}^l(s, a) \big)_{(s, a) \in \mathcal{K}}$, respectively, are known.   
  \item \label{mean_UB_
LB_RC_Constr_Known} 
The vector $\mu_{\tilde{d}^k}$, and upper and lower bound vectors of $\tilde{d}^k$, denoted by 
$\tilde{d}^{ku} = \big( \tilde{d}^{ku}(s, a) \allowbreak \big)_{(s, a) \in \mathcal{K}}$ and $\tilde{d}^{kl}= \big( \tilde{d}^{kl}(s, a) \big)_{(s, a) \in \mathcal{K}}$, respectively, are known. 
  \item \label{Compo_RC_independent}  The components of $\tilde{c}$ and $\tilde{d}^k$, $k \in \mathbb{K}$, are independent.
    \item \label{sub_gaussian_RC} For each $(s,a) \in \mathcal{K}$, $\tilde{c}(s,a)$ and  $\tilde{d}^k(s,a)$ are  sub-Gaussian random variables with mean vectors $\mu_{\tilde{c}}(s,a)$ and   $\mu_{\tilde{d}^k}(s,a)  $,  and  have sub-Gaussian parameters $\sigma_{\tilde{c}} (s,a)$ and   $\sigma_{\tilde{d}^k} (s,a)$, respectively.
    \end{enumerate}
\end{assumption}

\subsection{Upper bound approximations}\label{Inner approximations for random costs}
We first introduce two lemmas that are used in subsequent analysis, then derive four different upper bound approximations, which are applicable under different conditions.

\begin{lemma}\label{convexity_of_sqrt_function}
   The function
     $   \hat{f}(y) = \displaystyle \sqrt{ \frac{ p_1^{y^{\frac{1}{\theta}}}}{1 -  p_1^{y^{\frac{1}{\theta}}} }  },$
    is a convex function of $y \in [0,1]$. 
\end{lemma}
\begin{proof}
    Consider the function 
    \begin{align*}
        \tilde{f}(y) & = \ln \Bigg( \sqrt{ \frac{ p_1 ^{y^{\frac{1}{\theta}}}}{1 -  p_1^{y^{\frac{1}{\theta}}} }  } \Bigg) \\ 
        & = \frac{1}{2}\bigg( y^{\frac{1}{\theta}} \ln p_1 - \ln \big( 1 -  p_1^{y^{\frac{1}{\theta}}} \big) \bigg)  = \frac{1}{2} \big( \tilde{f}_1(y) + \tilde{f}_2(y) \big). 
    \end{align*}
    The functions $\tilde{f}_1(y)$ and $p_1^{y^{\frac{1}{\theta}}}$  are convex in $y$, and the function $- \ln (1-x)$ is convex in $x$ because their second-order derivatives are non-negative.
     In addition,   $- \ln (1-x)$  is increasing in $x$ which makes  $\tilde{f}_2(y)  $ a convex function in $y$. Thus, $\tilde{f}(y)$ is a convex function of $y$, which in turn implies that  $\hat{f}(y)$ is convex in $y$ because  $ e^x$ is convex and 
increasing in $x$.
\end{proof}
\begin{lemma}\label{convexity_of_sqrt_function_hoeffding_term}\textnormal{(\!\cite{peng2022bounds})}
    Let $y \in [0,1]$ and $p_1 \geq 1-e^{-\frac{1}{2}} \approx 0.3935$. Then, the   function 
     $   \bar{f}(y) = \sqrt{- \frac{1}{2} \ln \big( 1 -  p_1 ^{y^{\frac{1}{\theta}}} \big) }$ is    
     a convex function of $y$. 
\end{lemma}
\begin{proof}
    Define $\displaystyle \bar{f}_1(x) = \sqrt{- \frac{1}{2} \ln \big( 1 -   x \big) }$, $x \in [p_1,1]$. For $x \in (p_1,1)$, its first and second-order derivatives are given by 
    \begin{align*}
        \bar{f}_1' (x) & = \frac{1}{4(1-x)} \bigg( -\frac{1}{2} \ln (1-x)  \bigg)^{-\frac{1}{2}} >0, \\
        \bar{f}_1'' (x) & = \frac{1}{8(1-x)^2} \bigg( -\frac{1}{2} \ln (1-x)  \bigg)^{-\frac{3}{2}} \bigg( -\ln(1-x) - \frac{1}{2} \bigg) \geq 0, 
    \end{align*}
    where the last inequality follows from the fact that $x > p_1 \geq 1-e^{-\frac{1}{2}}$. Hence, $\bar{f}_1$ is a convex and increasing function of $x$. 
    The function $ p_1 ^{y^{\frac{1}{\theta}}} $ is convex in $y$, thus the convexity of  $\bar{f}$ follows. 
\end{proof}
\begin{theorem}\label{C_H_Sub_G_random_costs_convex_approximation_statement}
   Let $f_0:(0,1) \rightarrow \mathbb{R}, V_0: [0,1] \rightarrow \mathbb{R}$, $f_k: (0,1)\times [0,1] \rightarrow \mathbb{R}$, and $V_k^*$ be a scalar, for all $k \in \mathbb{K}$. Consider the following general problem \begin{subequations}\label{C_H_Sub_G_random_costs_convex_approximation}
    \begin{align}
        \min_{z, \rho, (y_k)_{k \in \mathbb{K} } } & \ z \nonumber \\
   \textnormal{s.t.}     & \ \rho^T \mu_{\tilde{c}} + f_0(p_0)  V_0(\rho) \leq z,  \label{C_H_Sub_G_random_costs_objective}  \\ 
 & \ \rho^T \mu_{\tilde{d}^k} + f_k(p_1,y_k)  V_k^* \leq \xi_k, \ \forall \ k \in \mathbb{K}, \label{C_H_Sub_G_random_costs_constraints_convexified} \\ 
 & \
\sum_{k \in \mathbb{K} } y_k = 1, \  y_k \geq 0,    \ 
\rho \in \mathcal{Q}^{\alpha}(\gamma), \ \forall \ k \in \mathbb{K}. 
    \end{align}
    \end{subequations}
    
  \begin{enumerate}[wide]
 \item \label{One-sided Chebyshev inequality random costs convex approximation}
 Let  condition \ref{Mean_Cov_RC_known} of \Cref{assumption_costs_bounds_mean_known_independence}
holds. An upper bound to the optimal value of \eqref{ICCP_random_costs} is given by the optimal value of the convex programming problem \eqref{C_H_Sub_G_random_costs_convex_approximation}, when
$ \displaystyle f_0(p_0) = \sqrt{\frac{p_0}{1 - p_0}}$, $V_0(\rho) = \| \Sigma_{\tilde{c}}^{\frac{1}{2}} \rho \|_2$, $\displaystyle f_k(p_1, y_k) = \sqrt{ \frac{ p_1 ^{y_k^{\frac{1}{\theta}}}}{1 -  p_1 ^{y_k^{\frac{1}{\theta}}} }  }$, and $V_k^* =
 \displaystyle \max_{\bar{\rho} \in \mathcal{Q}^{ \alpha }(\gamma)} \sum_{ (s,a) \in \mathcal{K} } \allowbreak \bar{\rho} (
s,a)  \big\| \big(\Sigma_{\tilde{d}^k}^{\frac{1}{2}} \big)_{(s,a)}  \big\|_2$;  $\big(\Sigma_{\tilde{d}^k}^{\frac{1}{2}} \big)_{(s,a)}$ denotes the $(s,a)^{th}$ column vector of the matrix $ \Sigma_{\tilde{d}^k}^{\frac{1}{2}}  $,   $k \in \mathbb{K}$.  
 \item \label{Hoeffding inequality random costs convex approximation} Let conditions \ref{mean_UB_
LB_RC_Obj_Known}-\ref{Compo_RC_independent} of \Cref{assumption_costs_bounds_mean_known_independence} hold and $p_1 \geq 1-e^{-\frac{1}{2}}$. An upper bound to the optimal value of  \eqref{ICCP_random_costs} is given by the optimal value of the   convex programming  problem \eqref{C_H_Sub_G_random_costs_convex_approximation}, when  $\displaystyle f_0(p_0) = \sqrt{-\frac{1}{2}    
   \ln (1- p_0)} $, $V_0(\rho) = \| D_{\tilde{c}}^{\frac{1}{2}} \rho \|_2$, $ 
\displaystyle f_k(p_1, y_k) = \sqrt{- \frac{1}{2} \ln \big( 1 -  p_1 ^{y_k^{\frac{1}{\theta}}} \big) }$, and $V_k^* = \displaystyle \max_{\bar{\rho} \in \mathcal{Q}^{ \alpha }(\gamma)} \sum_{ (s,a) \in \mathcal{K} }  \bar{\rho} (s,a)  \big( \tilde{d}^{ku} (s,a) - \tilde{d}^{kl} (s,a) \big)$; $D_{\tilde{c}}$ is a $\vert \mathcal{K} \vert \times \vert \mathcal{K} \vert $ diagonal matrix whose $(s,a)^{th}$ diagonal entry is defined as  $\big(\tilde{c}^u(s,a) - \tilde{c}^l(s,a)\big)^2$.
\item \label{sub_G_theorem_statement_random_costs} Let conditions \ref{Compo_RC_independent}-\ref{sub_gaussian_RC} of \Cref{assumption_costs_bounds_mean_known_independence} hold and  $p_1 \geq 1-e^{-\frac{1}{2}}$. An upper bound to the optimal value of \eqref{ICCP_random_costs} is given by the optimal value of the  convex programming  problem \eqref{C_H_Sub_G_random_costs_convex_approximation}, when $f_0(p_0) = \sqrt{-2   
   \ln (1- p_0)} $, $V_0(\rho)= \| \Sigma_{\tilde{c}, SG}^{\frac{1}{2}} \rho \|_2$, $f_k(p_1, y_k) = \sqrt{- 2 \ln \big( 1 -  p_1 ^{y_k^{\frac{1}{\theta}}} \big) }$, and $V_k^* = \displaystyle \max_{\bar{\rho} \in \mathcal{Q}^{ \alpha }(\gamma)} \sum_{ (s,a) \in \mathcal{K} }  \bar{\rho} (s,a)  \sigma_{\tilde{d}^k} (s,a)$;  $\Sigma_{\tilde{c}, SG}$   
   is a $\vert \mathcal{K} \vert \times \vert \mathcal{K} \vert $ diagonal matrix whose $(s,a)^{th}$ diagonal entry is defined as  $\sigma^2_{\tilde{c}} (s,a)$. 

    \end{enumerate}
 \end{theorem}
 \begin{proof}
 \begin{enumerate}[wide=0pt]
      \item It is sufficient to show that the feasible region of \eqref{C_H_Sub_G_random_costs_convex_approximation} is contained in the feasible region of \eqref{ICCP_random_costs}. Let  $\big(z, \rho, (y_k)_{k \in \mathbb{K}} \big)$ be a feasible   solution of \eqref{C_H_Sub_G_random_costs_convex_approximation}. Then,  
\begin{align*}
    V_k^* & \geq   \sum_{ (s,a) \in \mathcal{K} }  \rho (s,a)  \big\| \big(\Sigma_{\tilde{d}^k}^{\frac{1}{2}} \big)_{(s,a)}  \big\|_2 \\ 
    & \geq   \big\| \sum_{ (s,a) \in \mathcal{K} }  \rho (s,a)   \big(\Sigma_{\tilde{d}^k}^{\frac{1}{2}} \big)_{(s,a)}  \big\|_2  = \| \Sigma_{\tilde{d}^k}^{\frac{1}{2}} \rho \|_2, \ \forall \ k \in \mathbb{K}, 
\end{align*}
where the first and the second inequalities follow by the  definition of $V_k^*$  and by triangular inequality, respectively. Thus,  the vector satisfies 
\begin{align}
 & \ \rho^T \mu_{\tilde{d}^k} + \sqrt{ \frac{ p_1^{y_k^{\frac{1}{\theta}}}}{1 -  p_1 ^{y_k^{\frac{1}{\theta}}} }  } \ \| \Sigma_{\tilde{d}^k}^{\frac{1}{2}} \rho \|_2 \leq \xi_k, \ \forall \ k \in \mathbb{K}. \label{One-sided_Chebyshev_inequality_random_costs_constraints}
\end{align}
From  part \ref{One-sided Chebyshev inequality} of 
\cref{C_H_Sub_G_common_statement_general_CC}, if $\big(z, \rho, (y_k)_{k \in \mathbb{K}} \big)$ is a feasible solution of 
\eqref{C_H_Sub_G_random_costs_objective} and \eqref{One-sided_Chebyshev_inequality_random_costs_constraints}, then it is a feasible solution of 
 \Crefrange{ICCP_random_costs_objective}{ICCP_random_costs_constraints}.  The convexity of \eqref{C_H_Sub_G_random_costs_convex_approximation} follows from \Cref{convexity_of_sqrt_function} and the fact that $V_k^* \geq 0$,  $k \in \mathbb{K}$.
 \end{enumerate}
     \begin{enumerate}[wide,resume]
 \item Similar to the  first part, let  $\big(z, \rho, (y_k)_{k \in \mathbb{K}} \big)$ be a feasible   solution of \eqref{C_H_Sub_G_random_costs_convex_approximation}. Then,  
 \begin{align*}
    V_k^* & \geq 
    \sum_{ (s,a) \in \mathcal{K} }  \rho (s,a)  \big( \tilde{d}^{ku} (s,a) - \tilde{d}^{kl} (s,a) \big) \\
     & =   \sum_{ (s,a) \in \mathcal{K} }  \rho (s,a)  \big\| \big( D_{\tilde{d}^k}^{\frac{1}{2}} \big)_{(s,a)}  \big\|_2 \\ 
    & \geq   \big\| \sum_{ (s,a) \in \mathcal{K} }  \rho (s,a)   \big( D_{\tilde{d}^k}^{\frac{1}{2}} \big)_{(s,a)}  \big\|_2 = \| D_{\tilde{d}^k}^{\frac{1}{2}} \rho \|_2, \ \forall \ k \in \mathbb{K}, 
\end{align*}
where   $D_{\tilde{d}^k }$  is a $\vert \mathcal{K} \vert \times \vert \mathcal{K} \vert $ diagonal matrix whose $(s,a)^{th}$ diagonal entry is defined as    $\big(\tilde{d}^{k u}(s,a) - \tilde{d}^{k l}(s,a)\big)^2$,   $k \in \mathbb{K}$. The remaining arguments are similar to the first part using part \ref{Hoeffding inequality} of \Cref{C_H_Sub_G_common_statement_general_CC}. The convexity of \eqref{C_H_Sub_G_random_costs_convex_approximation} follows from \Cref{convexity_of_sqrt_function_hoeffding_term} and the fact that $V_k^* \geq 0$,  $k \in \mathbb{K}$. 
 \item The result follows from arguments similar to the second part using part \ref{Sub-Gaussian condition} of \Cref{C_H_Sub_G_common_statement_general_CC}.
 \end{enumerate}
 \end{proof}
 \begin{remark}\label{comparison_f_k_values}
   The upper bound approximations under parts \ref{One-sided Chebyshev inequality random costs convex approximation}-\ref{sub_G_theorem_statement_random_costs} of \Cref{C_H_Sub_G_random_costs_convex_approximation_statement} differ in the values of $f_0(p_0)$, $V_0(\rho)$, $f_k(p_1,y_k)$, and $V_k^*$. For any $p\in (0,1)$ and $\theta\ge 1$ we can show that $ \displaystyle  \sqrt{- \frac{1}{2} \ln \big( 1 -  p^{w^{\frac{1}{\theta}}} \big) }  \leq   \displaystyle  \sqrt{ \frac{ p^{w^{\frac{1}{\theta}}}}{1 -  p^{w^{\frac{1}{\theta}}} }}$, for all $w\in[0,1]$. Additionally, when $p\ge  \frac{1}{2}$, $ \displaystyle  \sqrt{- 2 \ln \big( 1 - p^{w^{\frac{1}{\theta}}} \big) }  \leq   \displaystyle  \sqrt{ \frac{ p^{w^{\frac{1}{\theta}}}}{1 -  p^{w^{\frac{1}{\theta}}} }  }$,  for all  $w\in[0,1]$. Thus, if condition \ref{Compo_RC_independent} of \Cref{assumption_costs_bounds_mean_known_independence} holds and the sub-Gaussian parameters of costs $\tilde{c}(s,a)$ and $\tilde{d}^k(s,a)$ are same as their variances, the values of $V_k^*$ under parts \ref{One-sided Chebyshev inequality random costs convex approximation} and \ref{sub_G_theorem_statement_random_costs} are same. Consequently, the approximation under part \ref{sub_G_theorem_statement_random_costs} is better than part \ref{One-sided Chebyshev inequality random costs convex approximation}.    
\end{remark}
\begin{theorem}\label{Bernstein inequality random costs} 
    Let conditions \ref{mean_UB_
LB_RC_Obj_Known}-\ref{Compo_RC_independent} of \Cref{assumption_costs_bounds_mean_known_independence} hold. An upper bound to the optimal value of \eqref{ICCP_random_costs} is given by the optimal value of the following convex programming  problem 
        \begin{align}\label{Bernstein_inequality_random_costs}
         \min_{z, \rho, (y_k)_{k \in  \mathbb{K} } }  & z \nonumber \\
       \textnormal{s.t.}       
       \sum_{(s,a) \in \mathcal{K}}  &
       \begin{aligned}[t] 
       \ln \big\{  A_{\tilde{c}}(s,a)   \exp \big( h^{rc}_0 \rho(s,a)  \tilde{c}^u(s,a)  \big)  +   \big(  1 - A_{\tilde{c}}(s,a) \big)  \exp \big( h^{rc}_0  \rho(s,a)  \tilde{c}^l(s,a) \big) \big\} 
       &
       \\   \leq \ln (1-p_0) + h^{rc}_0 z,
       &
       \end{aligned} 
       \nonumber \\
        \sum_{(s,a) \in \mathcal{K}} &  \begin{aligned}[t] 
        \ln  \big\{  A_{\tilde{d}^k}(s,a)   \exp \big( h^{rc}_k  \rho(s,  a) \tilde{d}^{k u}  (s,a)  \big)  +   \big(  1 - A_{\tilde{d}^k}(s,a)  \big) \exp \big( h^{rc}_k   \rho(s, a)  \tilde{d}^{k l}(s, 
        &
        \\  a)  \big)   \big\}   \leq 
        \ln \big( 1 -  p_1 ^{y_k^{\frac{1}{\theta}}} \big)  
        + h^{rc}_k \xi_k, \ \forall \ k \in \mathbb{K},
        &
        \end{aligned} \nonumber \\ 
      & 
      \sum_{k \in \mathbb{K} }  y_k = 1, \ y_k \geq 0,  \  \rho \in \mathcal{Q}^{\alpha}(\gamma), \ \forall \ k \in \mathbb{K},
    \end{align}
    for  arbitrarily fixed constants $h^{rc}_0 > 0$, $h^{rc}_k >0$,   $k \in \mathbb{K}$, and where $\displaystyle A_{\tilde{c}}(s,a) = \frac{\mu_{\tilde{c}}(s,a) - \tilde{c}^l(s,a)  }{ \tilde{c}^u(s,a) - \tilde{c}^l(s,a) }  $, $\displaystyle A_{\tilde{d}^k}(s,a) = \frac{\mu_{\tilde{d}^k }(s,a) - \tilde{d}^{k l}(s,a)  }{ \tilde{d}^{k u}(s,a) - \tilde{d}^{k l}(s,a) } $,   $(s,a ) \in \mathcal{K}$,    $k \in \mathbb{K}$. 
\end{theorem}
\begin{proof}
    It is sufficient to show that the feasible region of \eqref{Bernstein_inequality_random_costs} is contained in the feasible region of \eqref{ICCP_random_costs}. This claim directly follows from \Cref{Bernstein inequality}. Furthermore, from the proof of \Cref{convexity_of_sqrt_function},   $\ln \big( 1 -  p_1 ^{y_k^{\frac{1}{\theta}}} \big)$ is a concave function in $y_k$. Since all other terms and constraints are convex,  the convexity of \eqref{Bernstein_inequality_random_costs} follows. 
\end{proof}

\subsection{Lower bound approximation}\label{Outer approximation for random costs}
To derive a lower bound approximation  of the JCCMDP problem, we use the following lemma.  
\begin{lemma}\label{taylor_approximation_costs_case_statement}
    For a fixed number $N$, let $\mathcal{N}$ = $\{ 1,2,\ldots, N \}$ and  $(y^i)_{i \in \mathcal{N} }$  be chosen from the interval $[0,1] $ such that $y^1 < y^2 < \ldots <y^N$. 
   Then, for $y \in [0,1]$, 
    \begin{align*}
        p_1^{y^{ \frac{1}{\theta} }} \geq \max_{i \in \mathcal{N}} \, \big( a^i + b^iy \big),
    \end{align*}
    where  
    $\displaystyle a^i = p_1^{ {y^i}^{\frac{1}{\theta}} } \bigg( 1 -  \frac{{y^i}^{\frac{1}{\theta}}}{\theta} \ln p_1 \bigg)$ and   $b^i = \displaystyle p_1^{ {y^i}^{\frac{1}{\theta}} }  \frac{{y^i}^{\frac{1}{\theta} - 1}}{\theta} \ln p_1  , \ \forall \ i \in \mathcal{N}$. 
\end{lemma}
\begin{proof}
    The first-order Taylor series expansion  of $p_1 ^{y^{ \frac{1}{\theta} }}$     
    around the point $y^i$, denoted by $( p_1 )^{y^{\frac{1}{\theta}}}_i$, is defined as 
    \begin{align*} 
   (p_1)^{y^{\frac{1}{\theta}}}_i
   & =
    p_1^{{y^i}^{\frac{1}{\theta}}} + (y - y^i) \frac{d}{dy} p_1^{{y^i}^{\frac{1}{\theta}}},   \\ 
   & =
   p_1^{ {y^i}^{\frac{1}{\theta}} } \bigg( 1 -  \frac{{y^i}^{\frac{1}{\theta}}}{\theta} \ln p_1 \bigg) + 
    y p_1^{ {y^i}^{\frac{1}{\theta}} }  \frac{{y^i}^{\frac{1}{\theta} - 1}}{\theta}  \ln p_1    =
   a^i + b^i y, \ \forall \ i \in \mathcal{N}. 
\end{align*}
 Since $ p_1 ^{y^{\frac{1}{\theta}}}$ is convex in $y$, $ \displaystyle  p_1 ^{y^{\frac{1}{\theta}}} \geq \max_{i \in \mathcal{N}} \, (p_1)^{y^{\frac{1}{\theta}}}_i$ and the result follows. 
\end{proof}
\begin{theorem}\label{linear_approximation_random_costs_statement}
Let conditions \ref{mean_UB_
LB_RC_Obj_Known}-\ref{mean_UB_
LB_RC_Constr_Known} of \Cref{assumption_costs_bounds_mean_known_independence} hold. A lower bound  to the optimal value of \eqref{ICCP_random_costs} is given by the optimal value of the following LP problem 
\begin{subequations}\label{linear_approximation_random_costs}
        \begin{align}
            & \min_{ z, \rho, (\bar{y}_k)_{k \in \mathbb{K} }, \mathfrak{m}^{rc}_{c}, \mathfrak{m}^{rc}_{d} }  \ z \nonumber \\
              \textnormal{s.t.} 
             & \ 
             \rho^T \mu_{\tilde{c}} - z \leq (1-p_0) \mathfrak{m}^{rc}_c, \ 
              \lambda^{rc}_c \leq \mathfrak{m}^{rc}_c,  \ 
             \rho^T \tilde{c}^u - z \leq \mathfrak{m}^{rc}_c,    \ 
            \rho^T \tilde{c}^l \leq z, \label{linear_approximation_random_costs_objective_using_expectation_quantile}  \\ 
            & \
            \rho^T \mu_{\tilde{d}^k } -\xi_k + a_k^i \mathfrak{m}^{rc}_d + b_k^i \bar{y}_k \leq \mathfrak{m}^{rc}_d,  \ \forall \ i \in \mathcal{N},  \ k \in \mathbb{K}, \label{linear_approximation_random_costs_constraints_taylor_approximation}\\ 
           & \  
            \lambda^{rc}_d \leq \mathfrak{m}^{rc}_d, 
            \  \rho^T \tilde{d}^{ku} - \xi_k \leq \mathfrak{m}^{rc}_d, \ 
           \rho^T \tilde{d}^{kl} \leq \xi_k, \ \forall \ k \in \mathbb{K}, \label{linear_approximation_random_costs_constraints_quantile} \\ 
            & \
        \sum_{k \in \mathbb{K} } \bar{y}_k = \mathfrak{m}^{rc}_d, \  \bar{y}_k \geq 0,     \ 
        \rho \in \mathcal{Q}^{\alpha}(\gamma), \ \forall \ k \in \mathbb{K},  \label{linear_approximation_random_costs_sum_1_non_negative_flow_balance}
        \end{align}
    \end{subequations}
    where $\lambda^{rc}_c > 0$, $\lambda^{rc}_d >0$,  are constants and 
    \begin{align}\label{definition_aki_bki}
    \displaystyle  a_k^i  = p_1^{ {y_k^i}^{\frac{1}{\theta}} } \bigg( 1 -  \frac{{y_k^i}^{\frac{1}{\theta}}}{\theta} \ln p_1 \bigg),   \ 
      b_k^i  =  p_1^{ {y_k^i}^{\frac{1}{\theta}} }  \frac{{y_k^i}^{\frac{1}{\theta} - 1}}{\theta} \ln p_1,  \ \forall \ i \in \mathcal{N},  \ k \in \mathbb{K}.
    \end{align}
\end{theorem}
\begin{proof}
It is sufficient to show that for every feasible solution of \eqref{ICCP_random_costs}, there exists a feasible solution of \eqref{linear_approximation_random_costs}. 
 Let $\big(z, \rho, (y_k)_{k \in \mathbb{K}} \big)$ be a feasible   solution of \eqref{ICCP_random_costs}.  
    Define  $\mathfrak{m}^{rc}_c = \max \big(\lambda^{rc}_c, \rho^T \tilde{c}^u - z \big)$ and $ \displaystyle \mathfrak{m}^{rc}_d =  \VV{\max} \big(\lambda^{rc}_d, \VV{(\rho^T \tilde{d}^{ku} - \xi_k)_{k \in \mathbb{K}}} \big)$. From \Cref{general lower bound}, \eqref{linear_approximation_random_costs_objective_using_expectation_quantile} and  \eqref{linear_approximation_random_costs_constraints_quantile} are satisfied  along with the following constraint
\begin{align}\label{outer_approximation_random_costs_constraint_product_terms}
    \rho^T \mu_{\tilde{d}^k } -\xi_k 
            +  p_1 ^{y_k^{\frac{1}{\theta}}} \mathfrak{m}^{rc}_d
            \leq  \mathfrak{m}^{rc}_d, \ \forall \ k \in \mathbb{K}.  
\end{align}
  For each $k \in \mathbb{K}$, let $y_k^1 < y_k^2 < \ldots < y_k^N$ be $N$ tangent points from $[0,1]$. It follows from \Cref{taylor_approximation_costs_case_statement} that 
 \begin{align*}
      p_1 ^{y_k^{\frac{1}{\theta}}}  \geq \max_{i \in \mathcal{N}} \, \big( a_k^i + b_k^i y_k \big),
 \end{align*}
 where $a_k^i$ and $b_k^i$ are defined as in \eqref{definition_aki_bki}. 
 Since $\big( \rho, (y_k)_{k \in \mathbb{K}} \big)$ satisfies  \eqref{outer_approximation_random_costs_constraint_product_terms},  it also  satisfies  
\begin{align*}
    \rho^T \mu_{\tilde{d}^k } -\xi_k 
            + (a_k^i + b_k^i y_k) \mathfrak{m}^{rc}_d
            \leq  \mathfrak{m}^{rc}_d, \ \forall  \ i \in \mathcal{N}, \ k \in \mathbb{K} . 
\end{align*}
Defining $\bar{y}_k = y_k \mathfrak{m}^{rc}_d$,  $k \in \mathbb{K}$, we observe that $\big(  \rho, (\bar{y}_k)_{k \in \mathbb{K}} \big)$ satisfies  \eqref{linear_approximation_random_costs_constraints_taylor_approximation} and \eqref{linear_approximation_random_costs_sum_1_non_negative_flow_balance}. Hence, 
  $\big( z, \rho, (\bar{y}_k)_{k \in \mathbb{K} }, \mathfrak{m}^{rc}_c, \mathfrak{m}^{rc}_d \big)$ is a feasible  solution of \eqref{linear_approximation_random_costs}. 
\end{proof}
\VV{
\begin{remark}
    In terms of state-action pairs, the time required to solve the above LP problem can be bounded above by $O( \vert \mathcal{K} \vert^{\frac{1}{2}} \vert \ln \hat{\epsilon} \vert )$, where $\hat{\epsilon}>0$ denotes the threshold below which the duality gap is to be reduced \cite{wright1997primal}. Assuming same action \VS{set $A$} at each state, we obtain $O( \vert S \vert^{\frac{1}{2}} \vert A \vert^{\frac{1}{2}} \vert  \ln \hat{\epsilon} \vert )$  \VS{as} the   \VS{upper bound} on time. 
\end{remark}
}
\subsection{Bounds on gap}\label{Theoretical guarantees costs}
We define the gap between the upper and lower bounds of the optimal value of the JCCMDP problem as 
 \begin{align}\label{gap_per}
    \text{Gap}(\%) =  \displaystyle \frac{ \text{ UB - LB }}{ \text{LB}}\times  100,
\end{align}
where UB and LB denote the optimal values of upper and lower bound approximations, respectively.
This gap is useful in determining the quality of approximations derived in Sections \ref{Inner approximations for random costs} and \ref{Outer approximation for random costs}. A lower value of \text{Gap}(\%)  indicates a better quality of approximations. The theoretical bounds on the gap can be derived using the extremal bounds on the optimal values of the upper and lower bound approximations. Let UB$_{rc}^{(u)}$ denote the upper bound on the optimal value of upper bound approximation and  LB$_{rc}^{(l)}$ denote the lower bound on the optimal value of lower bound approximation. Then,  \text{Gap}(\%)$\in [0, G_{rc}]$, where  $G_{rc} = \displaystyle \frac{ \text{UB}_{rc}^{(u)} \text{ - } \text{LB}_{rc}^{(l)} }{ \text{LB}_{rc}^{(l)}}\times  100$.
 Assuming that the approximations are feasible, we summarize the values of UB$_{rc}^{(u)}$ and LB$_{rc}^{(l)}$, in \Cref{Theoretical_Guarantees_costs}, corresponding to the
 approximations from Sections \ref{Inner approximations for random costs} and \ref{Outer approximation for random costs}. The extremal bounds UB$_{rc}^{(u)}$ and LB$_{rc}^{(l)}$ are derived by exploiting the structure of the respective problems. For example, the objective function value of the upper bound approximation corresponding to part \ref{One-sided Chebyshev inequality random costs convex approximation} of Theorem \ref{C_H_Sub_G_random_costs_convex_approximation_statement} is given by $\rho^T \mu_{\tilde{c}} + \sqrt{\frac{p_0}{1-p_0}}  \| \Sigma_{\tilde{c}}^{\frac{1}{2}} \rho \|_2 $ and it can be upper bounded by the corresponding term given in Table \ref{Theoretical_Guarantees_costs}. We use similar arguments to derive extremal bounds in all the cases. 
 \begingroup
    \small
\renewcommand\arraystretch{1.6}
    \begin{longtable}{c p{10.9cm}}
    \caption{Extremal bounds.}\\
    \endfirsthead
\caption{Extremal bounds.}\\
\midrule
\endhead
\midrule
\multicolumn{2}{r}{\footnotesize\itshape Continued on the next page}
\endfoot
\endlastfoot
\toprule
    Upper bounds & \multicolumn{1}{c}{UB$_{rc}^{(u)}$} \\ 
    \midrule
    \eqref{C_H_Sub_G_random_costs_convex_approximation} & \\ 
    \cline{1-1}
 Part \ref{One-sided Chebyshev inequality random costs convex approximation} 
    & 
    $\displaystyle \max_{\rho \in \mathcal{Q}^{\alpha}(\gamma)} \rho^T \mu_{\tilde{c}} + \sqrt{\frac{p_0}{1 - p_0}}  
    \bigg(
    \sum_{ (s,a) \in \mathcal{K} }  \rho (s,a)  \big\| \big(\Sigma_{\tilde{c}}^{\frac{1}{2}} \big)_{(s,a)}  \big\|_2
    \bigg)$
    \\
    Part \ref{Hoeffding inequality random costs convex approximation} 
    & 
    $\displaystyle \max_{\rho \in \mathcal{Q}^{\alpha}(\gamma)} \rho^T \mu_{\tilde{c}} + \sqrt{-\frac{1}{2}    
   \ln (1- p_0)} \bigg( \sum_{ (s,a) \in \mathcal{K} }  {\rho} (s,a)  \big( \tilde{c}^{u} (s,a) - \tilde{c}^{l} (s,a) \big) \bigg)$
    \\
    Part \ref{sub_G_theorem_statement_random_costs} 
    & 
    $\displaystyle \max_{\rho \in \mathcal{Q}^{\alpha}(\gamma)} \rho^T \mu_{\tilde{c}} + \sqrt{-2   
   \ln (1- p_0)} \bigg( \sum_{ (s,a) \in \mathcal{K} }   {\rho} (s,a)  \sigma_{\tilde{c}} (s,a) \bigg)$\\
   \cline{1-1}
   \eqref{Bernstein_inequality_random_costs} &
   $\displaystyle \max_{\rho \in \mathcal{Q}^{\alpha}(\gamma)}  \rho^T \tilde{c}^u -  \frac{1}{h^{rc}_0} \ln (1-p_0)$ 
   \\
   \midrule
   Lower bound & \multicolumn{1}{c}{LB$_{rc}^{(l)}$} \\ 
    \midrule
    \eqref{linear_approximation_random_costs} & $\displaystyle \min_{\rho \in \mathcal{Q}^{\alpha}(\gamma)}  \rho^T \tilde{c}^l$
   \\
   \bottomrule
   \label{Theoretical_Guarantees_costs}
    \end{longtable}
    \endgroup
Using these extremal bounds, we can calculate $G_{rc}$, which in turn gives an interval for Gap(\%). However, the extremal bounds need not be precise because they are derived using a subset of the constraints from the feasible region of the respective problem. In the numerical experiments, we observe a significant reduction in the gap  obtained
by solving our approximations as compared to the theoretical ones.

\section{JCCMDP problem under random running costs and  transition probabilities}\label{JCCMDP problem under random running costs and transition probabilities}
In this section, we consider the JCCMDP problem for the case \ref{JCCMDP_transition_probabilities_desc.}, i.e.,  both running costs and the transition probabilities are random.  In Lemma \ref{lemma}, we show that the JCCMDP problem for the case \ref{JCCMDP_running_costs_desc.} can be restricted to the class of stationary policies without loss of optimality. This is true because, for fixed transition probabilities, there is a one-to-one correspondence between the set of stationary policies and the set $\mathcal{Q}^\alpha(\gamma)$ defined by \eqref{definition_q_alpha}. When the transition probabilities are also random,  this approach cannot be used.
The presence of random transition probabilities within a chance constraint makes the problem difficult because
the expected discounted costs, even for stationary policies, are non-linear (in fact, non-convex) functions of the transition probabilities. The presence of history-dependent policies adds another level of difficulty. Therefore, in this paper, we study the JCCMDP problem \eqref{common_JCCMDP} when $\mathfrak{v} = rc{\text-}tp$  restricted to the class of stationary policies\VV{, defined as follows:
\begin{subequations}\label{common_JCCMDP_uncertain_TP}
\begin{align}
         & \min_{z,f \in F_{S}} \  z \nonumber \\
        \text{s.t. } & \mathbb{P} \Big(  \tilde{C}_{rc{\text-}tp} \big( \gamma,f \big)  \leq z \Big)\geq p_0, \label{common_JCCMDP_uncertain_TP_objective} \\
         & \mathbb{P}\Big(   \tilde{D}_{rc{\text-}tp}^k \big( \gamma,f \big)  \leq \xi_k,  \ \forall \ k \in \mathbb{K} \Big)\geq p_1. \label{JCC_common_JCCMDP_uncertain_TP}
\end{align}
\end{subequations}} 
\noindent
For a given  $f \in F_S$, the random expected discounted costs can be defined in matrix form as  \cite{wiesemann2013robust}
\begin{align}\label{vector_expression_costs}
 \tilde{C}_{rc{\text-}tp} \big( \gamma,f \big )
 & =
 (1-\alpha) \gamma^T   \tilde{Q}_f \tilde{c}_f, \nonumber\\
 \tilde{D}_{rc{\text-}tp}^k \big( \gamma,f \big ) 
 & =
 (1-\alpha) \gamma^T \tilde{Q}_f \tilde{d}^k_f, \ \forall \ k \in \mathbb{K}, 
\end{align}
where $\tilde{c}_f \in \mathbb{R}^{\vert S \vert }$ and $\tilde{d}^k_f \in \mathbb{R}^{\vert S \vert } $ denote the   random running cost vectors for a given policy $f$ whose $s^{th}$ components are defined as  $\tilde{c}_f(s) =  \displaystyle \sum_{a \in A(s)} f(s, a) \tilde{c}( s, a) $ and $\tilde{d}^k_f(s) = \displaystyle \sum_{a \in A(s)} f(s, \allowbreak a) \tilde{d}^k( s, \allowbreak a)$,  $k \in \mathbb{K}$, respectively.  Furthermore,  $\tilde{Q}_f(\omega) = \big(I -  \alpha \tilde{P}_{ f }(\omega) \big)^{-1}$,   $\omega \in \Omega$, where $I \in \mathbb{R}^{\vert S \vert \times \vert S \vert} $ denotes the identity matrix  and  $\tilde{P}_f \in \mathbb{R}^{\vert S \vert \times \vert S \vert}$ denotes the random transition probability matrix induced by $f$. For a transition from a state  $s $ to  $s'$,
the component of $\tilde{P}_f (\omega)$ is given by $ \displaystyle \sum_{a \in A(s)} f(s, a) \tilde{p}(s' \vert s, a) (\omega) $.  
For each $(s,a) \in \mathcal{K}$, $s' \in S$, and $\omega \in \Omega$, we assume that we can express $\tilde{p}(s' \vert s, a)(\omega)$  as
\begin{align}\label{expression_random_TP_perturbation_matrix}
    \tilde{p}(s' \vert s, a)(\omega) = \mu(s' \vert s, a) + \zeta(s' \vert s, a)(\omega), 
\end{align}
where $\mu(s' \vert s, a)$ and $\zeta (s' \vert s, a)$  denote the expected value of $\tilde{p}(s' \vert s, a)$ and the random perturbation of $\tilde{p}(s' \vert s, a)$ from  $\mu(s' \vert s, a)$, respectively.
The random perturbation $\zeta (s' \vert s, a)$ is such that for every $\omega\in \Omega$,  both $\tilde{p}(s' \vert s,a)(\omega)$ and $\mu (s' \vert s,a)$, are transition probabilities, i.e., for all $(s,a) \in \mathcal{K},  s' \in  S$,
 $\tilde{p}(s' \vert s,a)(\omega)\ge 0$, $\mu (s' \vert s,a) \geq 0$, and $ \displaystyle \sum_{s' \in S} \tilde{p}(s' \vert s,a)(\omega)$ = $ \displaystyle \sum_{s' \in S} \mu (s' \vert s,a)$ = 1. 
 Thus, the random perturbations $ \big\{ \zeta (s' \vert s,a) \allowbreak \big\}_{ (s,a) \in \mathcal{K}, s' \in S } $ must be bounded such that  $\displaystyle \sum_{s'\in S} \zeta (s' \vert s,a)(\omega)=0$ and $\mathbb{E}( \zeta (s' \vert s,a) ) = 0$, for all   $(s,a) \in \mathcal{K}$ and  $s' \in S$. 
Let  $Z_f$ and $M_f$ denote the random perturbation matrix and the expected transition probability matrix, respectively.  Then using \eqref{expression_random_TP_perturbation_matrix}, $\tilde{P}_f$, $Z_f$, and $M_f$ are related as  $\tilde{P}_{ f }(\omega) = M_f + Z_f (\omega) $,  $\omega \in \Omega$, such that  for a transition from a state  $s $ to  $s'$,
the components of  $M_f$ and $Z_f(\omega)$ are defined as  $ \displaystyle \sum_{a \in A(s)} f(s, a) \mu  (s' \vert s, a)$ and $ \displaystyle \sum_{a \in A(s)} f(s, a) \zeta  (s' \vert s, a) (\omega)$, respectively.
Subsequently, we express the random expected discounted costs given in \eqref{vector_expression_costs} as a sum of two terms: the first term is the random expected discounted cost under random running costs and expected transition probability matrix $M_f$, while the second term involves the random perturbations.
\begin{proposition}\label{expression_inverse_m_matrix}
Let $ Q^M_f = \big( I - \alpha M_{ f } \big)^{-1} $. For a given  $f \in F_S$, 
\begin{align}\label{equivalent_expression_inverse_m_matrix}
\tilde{Q}_f(\omega)
=
Q^M_f \left( I 
+
\alpha  Z_{  f }(\omega) \tilde{Q}_f(\omega) \right), \ \forall \ \omega \in \Omega.
\end{align}
\end{proposition}
\begin{proof}
For a given $f \in F_S$, 
since $\tilde{P}_f(\omega)$ and $ M_{  f } $  are transition probability matrices,  $\tilde{Q}_f (\omega) $ and $Q^M_f$  exist for all $\omega \in \Omega$,  and we obtain  
\begin{align}\label{expression_for_existence_of_inverse}
\tilde{Q}_f (\omega)
& =  \left(  \left(I - \alpha Z_f (\omega) Q^M_f \right) (I - \alpha M_f) \right)^{-1} \nonumber  \\ 
& = Q^M_f
\left(I - \alpha Z_{  f } (\omega) Q^M_f \right)^{-1}, \ \forall \ \omega \in \Omega. 
\end{align}
  Furthermore, using the fact that  
$ \displaystyle 
\left( I - \alpha Z_{  f }(\omega) Q^M_f \right)
\left( I - \alpha Z_{  f }(\omega) Q^M_f \right)^{-1}
= I,
$
we obtain
\begin{align*}
 \left( I - \alpha Z_{  f }(\omega) Q^M_f \right)^{-1}
& = 
I + \alpha Z_{ f }(\omega) Q^M_f \left( I - \alpha Z_{ f }(\omega) Q^M_f \right)^{-1} \\ 
 & = 
 I + \alpha Z_{ f }(\omega) \tilde{Q}_f(\omega), \ \forall \ \omega \in \Omega,
\end{align*}
where the last equality follows from \eqref{expression_for_existence_of_inverse}. 
Pre-multiplying by $Q^M_f$ and using \eqref{expression_for_existence_of_inverse} gives  \eqref{equivalent_expression_inverse_m_matrix}.  
\end{proof}
We use Theorems \ref{C_H_Sub_G_common_statement_general_CC}, \ref{Bernstein inequality}, and \ref{general lower bound} to construct convex approximations whose optimal values provide upper and lower bounds to the optimal value of \eqref{common_JCCMDP_uncertain_TP}. 
The aforementioned theorems are only applicable for individual linear chance constraints, a characteristic which \eqref{common_JCCMDP_uncertain_TP} does not satisfy. Therefore, we first approximate \eqref{common_JCCMDP_uncertain_TP} into two individual linear chance-constrained programming (ILCCP) problems whose optimal values provide upper and lower bounds to its optimal value. 
We introduce some notations involving upper and lower bound vectors of $\tilde{c}$ and $\tilde{d}^k$, $k \in \mathbb{K}$. These notations are used in the subsequent analysis.
\subsubsection*
{Notations:}
\begin{enumerate}[label = (\roman*)]
\item $\displaystyle c_{\max} = \max_{(s,a) \in \mathcal{K}} \tilde{c}^u(s,a)$,   $\displaystyle d^k_{\max} = \max_{(s,a) \in \mathcal{K}} \tilde{d}^{ku}(s,a)$, $\forall \ k \in \mathbb{K} $,  
\item $\displaystyle c_{\min} = \min_{(s,a) \in \mathcal{K}} \tilde{c}^l(s,a)$, $\displaystyle d^k_{\min} = \min_{(s,a) \in \mathcal{K}} \tilde{d}^{kl}(s,a)$, $\forall \ k \in \mathbb{K} $,  
\item $\displaystyle C_{\max}(s)  = \max_{a \in A(s)} \tilde{c}^u(s,a) + \frac{\alpha}{1-\alpha} c_{\max}, \ \forall \ s \in S$, 
\item $\displaystyle C_{\min}(s)  = \min_{a \in A(s)} \tilde{c}^l(s,a) + \frac{\alpha}{1-\alpha} c_{\min},  \ \forall \ s \in S$,  
\item $\displaystyle D^k_{\max}(s)  = \max_{a \in A(s)} \tilde{d}^{ku}(s,a) + \frac{\alpha}{1-\alpha} d^k_{\max}, \ \forall \ s \in S, \ k \in \mathbb{K}$, 
\item $\displaystyle D^k_{\min}(s)  = \min_{a \in A(s)} \tilde{d}^{kl}(s,a) + \frac{\alpha}{1-\alpha} d^k_{\min},  \ \forall \ s \in S, \ k \in \mathbb{K}$. 
\end{enumerate}
To construct the ILCCP problems, we derive upper and lower bounds to the components of the vectors, $\tilde{Q}_f(\omega) \tilde{c}_f $ and $\tilde{Q}_f(\omega) \tilde{d}^k_f $,  whose $s^{th}$ components are denoted by $\big( \tilde{Q}_f (\omega) \tilde{c}_f \big)(s) $ and $  \big( \tilde{Q}_f (\omega) \tilde{d}^k_f \big)(s) $, respectively.
\begin{proposition}\label{bounds_for_inverse_m_matrix_with_costs}
For all  $\omega \in \Omega$, $f \in F_S$,   
\begin{enumerate}[wide] 
    \item   $  \big( \tilde{Q}_f (\omega) \tilde{c}_f \big)(s) \in \left[ C_{\min}(s), C_{\max}(s) \right], \ \forall \ s\in S$,
    \item   $ \big( \tilde{Q}_f (\omega) \tilde{d}^k_f \big)(s) \in \left[ D^k_{\min}(s), D^k_{\max}(s) \right], \ \forall \ s\in S, \ k\in\mathbb{K}.$ 
\end{enumerate}
\end{proposition}
\begin{proof}
\begin{enumerate}[wide=0pt]
    \item 
For a given  $\omega \in \Omega$ and $f \in F_S$, since $(I-\alpha \tilde{P}_f(\omega))\tilde{Q}_f(\omega) = I$, we obtain  
\begin{align*}
   \tilde{Q}_f(\omega) \tilde{c}_f  
   &  =
   \tilde{c}_f  + \alpha \tilde{P}_f (\omega) \tilde{Q}_f (\omega) \tilde{c}_f 
   \\
 & \geq
 \tilde{c}_f^l  + \frac{\alpha}{1-\alpha} \tilde{P}_f (\omega) c_{\min} \mathbbm{1}_{\vert S \vert}  
 =
 \tilde{c}_f^l  + \frac{\alpha}{1-\alpha} c_{\min} \mathbbm{1}_{\vert S \vert},  
\end{align*}
where $\mathbbm{1}_{\vert S \vert} \in \mathbb{R}^{ \vert S \vert  }$ denotes a  vector of ones. The first  inequality follows by the definition of $c_{\min}$ and  the fact that $\displaystyle \tilde{Q}_f (\omega) \mathbbm{1}_{\vert S \vert} = \frac{1}{1-\alpha}\mathbbm{1}_{\vert S \vert}$, while the last equality follows from the fact that $\displaystyle \tilde{P}_f (\omega)  \mathbbm{1}_{\vert S \vert} = \mathbbm{1}_{\vert S \vert}$.  Thus, 
\begin{align*}
\big( \tilde{Q}_f(\omega)  \tilde{c}_f  \big)(s)
& \geq 
\tilde{c}_f^l  (s) + \frac{\alpha}{1-\alpha} c_{\min} \\
& \geq
\min_{a \in A(s)} \tilde{c}^l(s,a) + \frac{\alpha}{1-\alpha} c_{\min} = C_{\min}(s), \ \forall \ s \in S. 
\end{align*} 
Similarly, 
\begin{align*}
   \big( \tilde{Q}_f(\omega)  \tilde{c}_f  \big)(s) & \leq C_{\max}(s), \ \forall \ s\in S.
\end{align*}
\end{enumerate}
\begin{enumerate}[wide,resume]
\item   The result follows from arguments similar to the first part. \qedhere
\end{enumerate}
\end{proof}
We use the results of  \Cref{bounds_for_inverse_m_matrix_with_costs} to derive bounds on the expected discounted costs.
As discussed previously, since $ \left\{ \zeta (s' \vert s, a) \right\}_{ (s, a) \in \mathcal{K}, s'  \in S } $ is bounded, we denote its upper and lower bounds by $ \left\{ \zeta^u (s' \vert s, a) \right\}_{ (s, a) \in \mathcal{K}, s'  \in S } $ and $ \left\{ \zeta^l (s' \vert s, a) \right\}_{ (s, a) \in \mathcal{K}, s'  \in S}$, respectively. In addition, let  $Z_f^u$ and $Z_f^l$ be matrices whose components denote the upper and lower bounds of the components of the matrix $Z_f(\omega)$, respectively,  for all $\omega \in \Omega$. Thus, the components of  $Z_f^u $ and $Z_f^l $ corresponding to a pair $(s,s')\in S\times S$ are defined as
$\displaystyle \sum_{a \in A(s)} f(s, a) \zeta^u(s' \vert s, a)  $ and $  \displaystyle \sum_{a \in A(s)} f(s, a) \zeta^l(s' \vert s, a)  $,  respectively.  
 We propose three upper bound and one lower bound approximations. Each approximation assumes a subset of conditions listed in Assumption \ref{existence_of_perturbations} below.
\begin{assumption}\label{existence_of_perturbations}
 \begin{enumerate}[wide=0pt]
 \item \label{Covariance_TP_known} The covariance matrix of the   perturbation vector $ \big( \zeta (s' \vert s,  a) \big)_{ (s, a) \in \mathcal{K}, s'  \in S } $ is known.   
 \end{enumerate}
 \begin{enumerate}[wide,resume]
 \item  \label{non_negativity_perturb_TP_assumption}  For each $(s,a)\in \mathcal{K}$, $s'\in S$, $ \zeta^u (s' \vert s,a)\ge 0$ and $\zeta^l (s' \vert s,a) \le 0$. 
\item \label{independence_TP} 
The set of $|\mathcal{K}|$ random vectors, $\left\{ \big(\zeta(s'|s,a) \big)_{s'\in S} \vert (s,a)\in \allowbreak \mathcal{K} \right\}$, is independent. 
\item   \label{costs_known_TP_assumption}  The vectors $\tilde{c}^u$, $\tilde{c}^l$,   $\tilde{d}^{ku}$, and $\tilde{d}^{kl}$,  $k \in \mathbb{K}$, are known.
\end{enumerate}
\end{assumption}

\subsection{Upper bound approximations}\label{Inner approximations random TPs}
\VV{Using \eqref{vector_expression_costs} and \Cref{expression_inverse_m_matrix}, we can write $\tilde{C}_{rc{\text-}tp} \left( \gamma,f \right)$ as 
\begin{align*}
     \tilde{C}_{rc{\text-}tp} \left( \gamma,f \right)
& = 
(1-\alpha) \gamma^T Q^M_f \left( \tilde{c}_f + 
\alpha  \left( Z_{ f } - Z_{  f }^l \right) \tilde{Q}_f \tilde{c}_f  + 
\alpha Z_{  f }^l \tilde{Q}_f  \tilde{c}_f \right). 
 \end{align*}
 }
 Under conditions \ref{non_negativity_perturb_TP_assumption} and \ref{costs_known_TP_assumption} of  \Cref{existence_of_perturbations}, 
 \VV{its} upper bound is given by 
\begin{align}\label{upper_bound_costs_under_random_TPs}
    \tilde{C}_{rc{\text-}tp} \left( \gamma,f \right)
& \leq 
(1-\alpha) \gamma^T Q^M_f \left( \tilde{c}_f^u  + 
\alpha  \left( Z_{  f } - Z_{  f }^l \right) C_{\max} + \alpha Z_{  f }^l C_{\min} \right), 
\end{align}
where the  inequality follows from \Cref{bounds_for_inverse_m_matrix_with_costs} along with the fact that  $Z_{  f }^l \leq \min\{0, Z_{  f }(\omega) \}   $,  $\omega \in \Omega$.
Similarly,  
\begin{align}\label{upper_bound_constraints_costs_under_random_TPs}
    \tilde{D}^k_{rc{\text-}tp} \big( \gamma,f \big)
& \leq 
(1-\alpha) \gamma^T Q^M_f \left( \tilde{d}^{ku}_f  + 
\alpha  \left( Z_{  f } - Z_{  f }^l \right) D^k_{\max} + \alpha Z_{  f }^l D^k_{\min} \right),  \ \forall \ k \in \mathbb{K}. 
\end{align}
\begin{lemma}\textnormal{(\!\cite{nemirovski2007convex})} \label{individual_constr_implies_joint_constr}
    If $f \in F_S$ is feasible for the following set of individual chance constraints 
    \begin{align} \label{ICC_implies_JCC}
   \mathbb{P}  \left( \tilde{D}_{rc{\text-}tp}^k \left( \gamma, f \right) \leq \xi_k \right)  \geq 1- \frac{1-p_1}{K}, \ \forall \ k \in \mathbb{K},  
\end{align}
then it is feasible for \eqref{JCC_common_JCCMDP_uncertain_TP}. 
\end{lemma}
\begin{proof} 
     Let $f \in F_S$ be feasible for  \eqref{ICC_implies_JCC}, then it satisfies 
    \begin{align*}
        1 - p_1
        & = \sum_{k \in \mathbb{K}} \frac{1 - p_1}{K} 
         \geq \sum_{k \in \mathbb{K}}  \mathbb{P}  \left( \tilde{D}_{rc{\text-}tp}^k \left( \gamma, f \right) > \xi_k \right) \\ 
        & \ \geq \mathbb{P}  \Big( \displaystyle \bigcup \limits_{k \in \mathbb{K}} \big\{ \tilde{D}_{rc{\text-}tp}^k \big( \gamma, f \big) > \xi_k \big\} \Big)  = 1 - \mathbb{P}  \left( \tilde{D}_{rc{\text-}tp}^k \left( \gamma, f \right) \leq \xi_k,  \ \forall \ k \in \mathbb{K} \right).
    \end{align*}
    Therefore, the result follows.  
\end{proof} 
Using \crefrange{upper_bound_costs_under_random_TPs}{upper_bound_constraints_costs_under_random_TPs} and \Cref{individual_constr_implies_joint_constr}, we construct an ILCCP problem whose optimal value is an upper bound to the optimal value of \eqref{common_JCCMDP_uncertain_TP}. 
\begin{lemma}\label{main-lemma}
Let conditions \ref{non_negativity_perturb_TP_assumption} and \ref{costs_known_TP_assumption} of  \Cref{existence_of_perturbations} hold.  
For every feasible solution of the ILCCP problem \eqref{ICC_inner_approximation_using_rho_removing_min_cost} given below, there exists a feasible solution of the JCCMDP problem  \eqref{common_JCCMDP_uncertain_TP}.  
\begin{subequations}\label{ICC_inner_approximation_using_rho_removing_min_cost}
\begin{align}
           \min_{ z, \, \rho } & \ z \nonumber \\
&          \textnormal{s.t.} 
 \ 
\mathbb{P}  \left( \sum_{(s,a) \in \mathcal{K}} \sum_{s' \in S \backslash s_c} \rho(s,a)   C_{\max}^{s_c}(s')  \zeta  (s' \vert s,a) \leq \mathcal{R}_c \right) \geq p_0, 
 \label{ICC_inner_approximation_objective} \\
&  
\begin{aligned}[t]
\mathbb{P}   \left( \sum_{(s,a) \in \mathcal{K}} \sum_{s' \in S \backslash s_{d^k}}  \rho(s,a)   D^{k, s_{d^k}}_{\max}(s')   \zeta  (s' \vert s,a) \leq \mathcal{R}_{d^k} \right)\geq  1 - \frac{1 - p_1}{K},  \forall  k \in \mathbb{K},
\end{aligned}
\label{ICC_inner_approximation_constraints} \\
& \ 
\rho \in \mathcal{Q}_M^{\alpha}(\gamma) \label{ICC_inner_approximation_flow_balance}, 
\end{align}
\end{subequations}
where $C_{\max}^{s_c}(s') = C_{\max}(s') - C_{\max}(s_c ) $, $D^{k, s_{d^k}}_{\max}(s') = D^k_{\max}(s') - D^k_{\max}(s_{d^k})$, for all $s' \in S$, $ \displaystyle s_c \in \argmin_{s' \in S} C_{\max}(s') $, $ \displaystyle s_{d^k} \in \argmin_{s' \in S} D^k_{\max}(s')$, and  $\mathcal{R}_c$, $\mathcal{R}_{d^k}$  are defined as 
\begin{align}\label{definition_r_c_and_r_d_k}
   & \alpha \mathcal{R}_c
     = 
    z -  \rho^T \tilde{c}^u  + \alpha \sum_{(s,a) \in \mathcal{K}} \sum_{s' \in S  }  \rho(s,a)  \left(  C_{\max}(s') - C_{\min}(s') \right) \zeta^l  (s' \vert s,a), \nonumber \\
    &\alpha \mathcal{R}_{d^k}
     =
    \xi_k -  \rho^T  \tilde{d}^{ku}  + \alpha \sum_{(s,a) \in \mathcal{K}} \sum_{s' \in S  }   \rho(s,a)   \left(  D^k_{\max}(s') - D^k_{\min}(s') \right) \zeta^l  (s' \vert s,a), \forall \ k \in \mathbb{K}, 
   \end{align} 
 and the set $\mathcal{Q}_M^{\alpha}(\gamma)$ is obtained from  \eqref{definition_q_alpha} by replacing $p({s}'|s,a)$ with $\mu({s}'|s,a)$,   $(s,a) \in \mathcal{K}$, $s' \in S$. 
\end{lemma}
\begin{proof}
    Let  $( z, \rho )$ be a feasible solution of \eqref{ICC_inner_approximation_using_rho_removing_min_cost}. 
Define 
\begin{equation}\label{policy}
     f(s,a) = \frac{\rho(s,a)}{\displaystyle \sum_{a \in A(s)} \rho(s,a)}, \ \forall \ (s,a)\in \mathcal{K},  
\end{equation} 
provided the denominator is non-zero. We arbitrarily choose $f(s)\in \wp(A(s))$   otherwise. Since $\rho \in \mathcal{Q}_M^{\alpha}(\gamma)$, the following equality holds \cite{altman}  
\begin{align}\label{rho_to_f}
 \rho(s,a) = \left( (1-\alpha) \gamma^T Q^M_f \right)(s) f(s,a), \ \forall \ (s,a) \in \mathcal{K}, 
\end{align}
where $\big( (1-\alpha) \gamma^T Q^M_f \big)(s)$ denotes the $s^{th}$-component of the vector $(1-\alpha) \gamma^T Q^M_f$. 
We obtain
\begin{align*}
     \alpha \mathcal{R}_c  & - \alpha \sum_{(s,a) \in \mathcal{K}} \sum_{s' \in S \backslash s_c} \rho(s,a)   C_{\max}^{s_c}(s')  \zeta  (s' \vert s,a) \\ 
   & = \alpha \mathcal{R}_c - \alpha \sum_{(s,a) \in \mathcal{K}} \sum_{s' \in S } \rho(s,a)   C_{\max}(s') \zeta  (s' \vert s,a)\\ 
& =  \begin{aligned}[t]
z -  \sum_{(s,a) \in \mathcal{K}}   \rho(s,a)  \bigg\{ \tilde{c}^u(s,a) + \alpha  \sum_{s' \in S }    \big( \zeta  (s' \vert s,a) -  \zeta^l  (s' \vert s,a) \big) C_{\max}(s')  \\ +\alpha \sum_{s' \in S } \zeta^l (s' \vert s,a)  C_{\min}(s')  \bigg\}
\end{aligned} \\ 
& = \begin{aligned}[t]
z -  \sum_{(s,a) \in \mathcal{K}} \big( (1-\alpha) \gamma^T Q^M_f \big)(s)  f(s,a)  \bigg\{ \tilde{c}^u(s,a) + \alpha  \sum_{s' \in S }    \big( \zeta  (s' \vert s,a) & \\  -\zeta^l  (s' \vert s,a) \big) C_{\max}(s') + \alpha \sum_{s' \in S } \zeta^l (s' \vert s,a)  C_{\min}(s')  \bigg\}&
\end{aligned} \\ 
   & =  z -  (1-\alpha) \gamma^T Q^M_f \left( \tilde{c}_f^u + 
\alpha  \left( Z_{  f } - Z_{  f }^l \right) C_{\max} + \alpha Z_{  f }^l C_{\min} \right)  \leq z - \tilde{C}_{rc{\text-}tp} \left( \gamma,f \right), 
\end{align*}
where the first equality follows using the fact that $\displaystyle \sum_{s' \in S} \zeta (s' \vert s, a)(\omega) = 0$,    $(s, a) \in \mathcal{K}$ and $\omega \in \Omega$, the second and the third equalities follow from \eqref{definition_r_c_and_r_d_k} and  \eqref{rho_to_f}, respectively. The last equality is the matrix formulation of its preceding expression, while the last inequality follows from   \eqref{upper_bound_costs_under_random_TPs}.  Therefore, $(z, f)$ is a feasible solution of \eqref{common_JCCMDP_uncertain_TP_objective}. Since $\rho$ is a feasible solution of \eqref{ICC_inner_approximation_constraints}, using arguments similar as above, $f$ defined by \eqref{policy} satisfies \eqref{ICC_implies_JCC}. Thus, from \Cref{individual_constr_implies_joint_constr}, $f$ satisfies \eqref{JCC_common_JCCMDP_uncertain_TP} and the result follows. 
\end{proof}
Since \eqref{ICC_inner_approximation_using_rho_removing_min_cost} is an ILCCP problem, we apply Theorems \ref{C_H_Sub_G_common_statement_general_CC} and  \ref{Bernstein inequality}   to obtain upper bounds to the optimal value of \eqref{common_JCCMDP_uncertain_TP}. 
\begin{theorem}\label{C_H_Sub_G_random_TPs_convex_approximation_thm_statement} 
Let $f_0, f_k: (0,1) \rightarrow \mathbb{R},  V_0, V_k: [0,1] \rightarrow \mathbb{R}$, for all $k \in \mathbb{K}$. Consider the following general problem
\begin{align}\label{C_H_Sub_G_random_TPs_convex_approximation}
     \min_{z, \, \rho} 
& \ z \nonumber \\
 \textnormal{s.t.}  & \ 
 f_0(p_0) \, V_0(\rho) \leq \mathcal{R}_c, \nonumber   \\
 & \
 f_k(p_1) \, V_k(\rho) \leq \mathcal{R}_{d^k}, \ \forall \  k \in \mathbb{K},  \
 \rho \in \mathcal{Q}_M^{\alpha}(\gamma).  
 \end{align}
     \begin{enumerate}[wide]
\item   \label{One-sided Chebyshev inequality random TPs convex approximation} 
Let conditions \ref{Covariance_TP_known}, \ref{non_negativity_perturb_TP_assumption}, and \ref{costs_known_TP_assumption} of \Cref{existence_of_perturbations}  hold.  An upper bound to the optimal value of \eqref{common_JCCMDP_uncertain_TP}
 is given by the optimal value of the   SOCP problem  \eqref{C_H_Sub_G_random_TPs_convex_approximation}, when $ \displaystyle f_0(p_0) = \sqrt{\frac{p_0}{1 - p_0}}$, $V_0(\rho) = \| \Sigma_{C \zeta} ^{\frac{1}{2}}  \rho \|_2$, $\displaystyle  f_k(p_1) = \sqrt{\frac{K-1+p_1}{1 - p_1}}$, and $V_k(\rho) = \| \Sigma_{D^k \zeta} ^{\frac{1}{2}}  \rho \|_2$; $ \Sigma_{C \zeta}$ and $\Sigma_{D^k \zeta}$ are covariance matrices of vectors whose $(s,a)^{th}$ components are given by $\displaystyle \sum_{s' \in S \backslash s_c}  C_{\max}^{s_c}(s') \zeta  (s' \vert s,a)$ and $ \displaystyle \sum_{s' \in S \backslash s_{d^k}}  D^{k, s_{d^k}}_{\max}(s')   \zeta (s' \vert s,  a)$, respectively. 
 \item   \label{Hoeffding inequality random TPs convex approximation}  Let conditions \ref{non_negativity_perturb_TP_assumption}-\ref{costs_known_TP_assumption} of \Cref{existence_of_perturbations} hold.  An upper bound to the optimal value of \eqref{common_JCCMDP_uncertain_TP}
 is given by the optimal value of the  SOCP problem \eqref{C_H_Sub_G_random_TPs_convex_approximation}, when $ 
\displaystyle f_0(p_0) = \sqrt{-\frac{1}{2} \ln   (1-p_0)  }$, $V_0(\rho) = \big\| \mathcal{D}_{C \zeta}^{\frac{1}{2}} \rho \big\|_2$, $\displaystyle f_k(p_1) = \sqrt{-\frac{1}{2} \ln  \frac{1 - p_1}{K } }$, $V_k(\rho) = \big\| \mathcal{D}_{ D^k \zeta}^{\frac{1}{2}} \rho \big\|_2$;  $\mathcal{D}_{C \zeta }  $, $\mathcal{D}_{D^k \zeta }  $   are  $\vert \mathcal{K} \vert \times \vert \mathcal{K} \vert $ diagonal matrices whose $(s,a)^{th}$ diagonal entries, denoted by,
$\big( \mathcal{D}_{C \zeta} \big)_{(s,a)}$ and $ \big( \mathcal{D}_{D^k \zeta} \big)_{(s,a)}$, respectively, are defined as
\begin{align*}
 & \big( \mathcal{D}_{C \zeta} \big)_{(s,a)}  = \Big( \displaystyle  \sum_{s' \in S \backslash s_c}  C_{\max}^{s_c}(s')  (  \zeta^u(s' \vert s,a) - \zeta^l (s' \vert s,a) )  \Big)^2, \, \forall \ (s,a) \in \mathcal{K}, \\ 
 & \begin{aligned}[t] \big( \mathcal{D}_{D^k \zeta} \big)_{(s,a)}  =
  \Big( \displaystyle  \sum_{s' \in S \backslash s_{d^k}}  D^{k, s_{d^k}}_{\max}(s')   (  \zeta^u(s' \vert s,a) - \zeta^l (s' \vert s,a) )\Big)^2,   \forall \ (s,a)   \in \mathcal{K}, 
  \ k \in \mathbb{K}.
  \end{aligned}
   \end{align*} 
     \end{enumerate}
\end{theorem}
\begin{proof}
    \begin{enumerate}[wide=0pt] 
    \item It is sufficient to show that the feasible region of \eqref{C_H_Sub_G_random_TPs_convex_approximation} is contained in the feasible region of \eqref{common_JCCMDP_uncertain_TP}. 
Let  $(z,\rho)$ be a feasible solution of  \eqref{C_H_Sub_G_random_TPs_convex_approximation}. The linear combinations of the random variables  $\big\{\big(\zeta(s'|s,a) \big)_{s'\in S} \mid (s,a) \allowbreak \in \mathcal{K}\big\}$, defined as,  
     $ \displaystyle \sum_{s' \in S \backslash s_c}  C_{\max}^{s_c}(s')  \zeta  (s' \vert s,a) $  and $\displaystyle \sum_{s' \in S \backslash s_{d^k}}  D^{k, s_{d^k}}_{\max}(s')  \zeta  ( \allowbreak s'  \vert s,a)$, for all  $(s,a) \in \mathcal{K}$, $k \in \mathbb{K}$, have expectations 0. It follows from  part \ref{One-sided Chebyshev inequality} of  \Cref{C_H_Sub_G_common_statement_general_CC} 
     that $(z, \rho)$ is a feasible solution of \eqref{ICC_inner_approximation_using_rho_removing_min_cost}. Then, from \Cref{main-lemma}, there exists a policy $f$ which is a feasible solution of \eqref{common_JCCMDP_uncertain_TP}.
     \end{enumerate}
     \begin{enumerate}[wide, resume]
     \item Let  $(z,\rho)$ be a feasible solution of  \eqref{C_H_Sub_G_random_TPs_convex_approximation}. It follows from  condition \ref{independence_TP} of \Cref{existence_of_perturbations} that the components of random vectors $ \allowbreak \displaystyle \Bigg(\sum_{s' \in S \backslash s_c}  C_{\max}^{s_c}(s') \zeta  (s' \vert s,a)\Bigg)_{(s,a)\in \mathcal{K}}$  and $\displaystyle \Bigg(\sum_{s' \in S \backslash s_{d^k}}     D^{k, s_{d^k}}_{\max}(s') \zeta  ( \allowbreak s'  \vert s,a)\Bigg)_{(s,a)\in \mathcal{K}}$,  
      $k \in \mathbb{K}$, are  independent. Moreover, the upper and lower bounds of their components  are $\displaystyle \sum_{s' \in S \backslash s_c}    C_{\max}^{s_c}(s') \zeta^u  (s' \vert s,a)$, $\displaystyle \sum_{s' \in S \backslash s_c}    C_{\max}^{s_c}( \allowbreak s')  \zeta^l  (s'  \vert s,a)$,  and $\allowbreak \displaystyle \sum_{s' \in S \backslash s_{d^k}}  \allowbreak  D^{k, s_{d^k}}_{\max}(s') \zeta^u  (  s' \vert s,a)$, $\allowbreak \displaystyle \sum_{s' \in S \backslash s_{d^k}}    D^{k, s_{d^k}}_{\max}(s') \zeta^l  (s' \vert s,a)$, respectively.
      It follows from  part \ref{Hoeffding inequality} of  \Cref{C_H_Sub_G_common_statement_general_CC} that $(z,\rho)$ is a feasible solution of \eqref{ICC_inner_approximation_using_rho_removing_min_cost}. Again, the result follows from \Cref{main-lemma}. \qedhere   
     \end{enumerate}
\end{proof}
\VV{
\begin{remark}
    In terms of state-action pairs, the time required to solve the SOCP problems given in the above theorem is bounded above by $O( \vert \mathcal{K} \vert^{\frac{7}{2}} )$ \cite{Lobo1998ApplicationsOS}. Assuming same action set $A$ at each state, \VS{the upper bound on time is given by} $O( \vert S \vert^{\frac{7}{2}} \vert A \vert^{\frac{7}{2}} )$. 
\end{remark}
}

\begin{theorem}\label{inner_approximation_using_bernstein_inequality_statement}
Let conditions \ref{non_negativity_perturb_TP_assumption}-\ref{costs_known_TP_assumption} of  \Cref{existence_of_perturbations} hold.  An upper bound to the optimal value of \eqref{common_JCCMDP_uncertain_TP}
 is given by the optimal value of the  following convex programming problem
\begin{align}\label{inner_approximation_using_bernstein_inequality_theorem_TP}
\min_{z, \, \rho} & \  z \nonumber  \\
   \textnormal{s.t.}  & 
\begin{aligned}[t]
   \sum_{(s,a) \in \mathcal{K}} 
  \ln \bigg\{    A_C (s,a)  \exp \Big( h_0^{rp}  \rho(s,a) \sum_{s' \in S \backslash s_c} C_{\max}^{s_c}(s')   \zeta^u  (s' \vert s,a) \Big) 
 + \big( 1 -  A_C (s,
 &
 \nonumber \\ 
 a) \big)  
   \exp \Big(  h_0^{rp}  \rho(s,a)\sum_{s' \in S \backslash s_c}  C_{\max}^{s_c}(s')  \zeta^l  (s' \vert s,a) \Big) \bigg\} 
  \leq   \ln  (1  -  p_0)  +  h_0^{rp} \mathcal{R}_c,
  &
  \end{aligned}
  \nonumber \\
   & 
   \begin{aligned}[b]
   \sum_{(s,a) \in \mathcal{K}}  
 \ln \bigg\{   A_{D^k} (s,a)   \exp \Big( h_k^{rp}  \rho(s,a) \sum_{s' \in S \backslash s_{d^k}}   D^{k, s_{d^k}}_{\max}(s')   \zeta^u (s' \vert s,  a) \Big) 
 +  
   \big( 1 - & \\
   A_{D^k} (s,a) \big)  \exp \Big(  h_k^{rp}  \rho(s,a) \sum_{s' \in S \backslash s_{d^k}}   D^{k, s_{d^k}}_{\max}(s') \zeta^l (s' \vert s,a) \Big) \bigg\} 
  \leq  
    \ln \frac{1 - p_1}{K} 
    &
    \\ 
    + h_k^{rp}  \mathcal{R}_{d^k}, \ \forall \ k \in \mathbb{K}, \  \rho  \in \mathcal{Q}_M^{\alpha}(\gamma),
    &
 \end{aligned}
\end{align}
 for some arbitrarily fixed $h^{rp}_0 > 0$, $h^{rp}_k >0$,  and where $\mathcal{R}_c$, $\mathcal{R}_{d^k}$,   $k \in \mathbb{K}$, are  defined as in \eqref{definition_r_c_and_r_d_k}, and
 \begin{align*}
 A_{C} ( s,a) & = 
    \frac{ - \displaystyle \sum_{s' \in S \backslash s_c} C_{\max}^{s_c}(s')  \zeta^l (s' \vert s,a)}{ \displaystyle \sum_{s' \in S \backslash s_c} C_{\max}^{s_c}(s')  \left( \zeta^u (s' \vert s,a) - \zeta^l (s' \vert s,a) \right) }, \ \forall \ (s,a) \in \mathcal{K}, \\ 
    A_{D^k} (s,a) & = \begin{aligned}[t]
    \frac{ - \displaystyle \sum_{s' \in S \backslash s_{d^k}} D^{k, s_{d^k}}_{\max}(s')  \zeta^l (s' \vert s,a)}{ \displaystyle \sum_{s' \in S \backslash s_{d^k}} D^{k, s_{d^k}}_{\max}(s') \left( \zeta^u (s' \vert s,a) - \zeta^l (s' \vert s,a) \right) }, \ \forall \ (s,a) \in \mathcal{K},  
    \ k \in \mathbb{K}. 
    \end{aligned}
    \end{align*}
\end{theorem}
\begin{proof}
The proof follows from \Cref{Bernstein inequality} and the arguments similar to the proof of 
part \ref{Hoeffding inequality random TPs convex approximation} of \Cref{C_H_Sub_G_random_TPs_convex_approximation_thm_statement}.    
 Since all the constraints are clearly convex, \eqref{inner_approximation_using_bernstein_inequality_theorem_TP} is a convex programming problem. 
\end{proof}

\subsection{Lower bound approximation}\label{Outer approximations random TPs}
Under conditions \ref{non_negativity_perturb_TP_assumption} and \ref{costs_known_TP_assumption} of  \Cref{existence_of_perturbations}, 
 a lower bound to $\tilde{C}_{rc{\text-}tp} \left( \gamma,f \right)$ is given by
 \begin{align}\label{lower_bound_costs_under_random_TPs}
    \tilde{C}_{rc{\text-}tp} \big( \gamma,f \big)
& = 
(1-\alpha) \gamma^T Q^M_f \left( \tilde{c}_f + 
\alpha  \left( Z_{ f } - Z_{  f }^l \right) \tilde{Q}_f \tilde{c}_f  + 
\alpha Z_{  f }^l \tilde{Q}_f  \tilde{c}_f \right) \nonumber \\ 
& \geq 
(1-\alpha) \gamma^T Q^M_f \left( \tilde{c}_f^l  + 
\alpha  \left( Z_{  f } - Z_{  f }^l \right) C_{\min} + \alpha Z_{  f }^l C_{\max} \right), 
\end{align}
where the  inequality follows from \Cref{bounds_for_inverse_m_matrix_with_costs} along with the fact that $Z_{  f }^l \leq \min\{0, Z_{  f }(\omega) \}   $,  $\omega \in \Omega$. Similarly, 
\begin{align}\label{lower_bound_constraints_costs_under_random_TPs}
      \tilde{D}^k_{rc{\text-}tp} \big( \gamma,f \big)
& \geq 
(1-\alpha) \gamma^T Q^M_f \left( \tilde{d}^{kl}_f  + 
\alpha  \left( Z_{  f } - Z_{  f }^l \right) D^k_{\min} + \alpha Z_{  f }^l D^k_{\max} \right), \ \forall\  k \in \mathbb{K}. 
\end{align}
Using \crefrange{lower_bound_costs_under_random_TPs}{lower_bound_constraints_costs_under_random_TPs}, we construct an ILCCP problem whose optimal value is a lower bound to the optimal value of \eqref{common_JCCMDP_uncertain_TP}. 
\begin{lemma}\label{lower-bd-lemma}
  Let conditions \ref{non_negativity_perturb_TP_assumption} and \ref{costs_known_TP_assumption} of  \Cref{existence_of_perturbations} hold.   For every feasible solution of the 
  JCCMDP problem \eqref{common_JCCMDP_uncertain_TP}, there exists a feasible solution of the following  ILCCP problem. 
\begin{subequations}\label{ICC_random_TPs_outer_approximation_simplified}
    \begin{align}
   &  \min_{ z, \, \rho }  \ z \nonumber \\
  &   \textnormal{s.t.}
      \  \mathbb{P}   \Bigg( \sum_{(s,a) \in \mathcal{K}} \sum_{s' \in S \backslash s_c' } \rho(s,a) C_{\min}^{s_c'}(s') \zeta  (s' \vert s,a) \leq \mathfrak{R}_c \Bigg) \geq p_0, \label{ICC_random_TPs_outer_approximation_objective_simplified}  \\ 
     & \
     \begin{aligned}[t]
      \mathbb{P}   \Bigg( \sum_{(s,a) \in \mathcal{K}} \sum_{s' \in S \backslash s_{d^k}' } \rho(s,a) D^{k, s_{d^k}'}_{\min}(s') \zeta  (s' \vert s,a)    \leq \mathfrak{R}_{d^k} \Bigg) \geq p_1,   \ \forall \ k \in \mathbb{K}, \label{ICC_random_TPs_outer_approximation_constraints_simplified} 
      \end{aligned}
      \\
     & \ \rho \in \mathcal{Q}_M^{\alpha}(\gamma) \label{ICC_random_TPs_outer_approximation_flow_balance_simplified},
     \end{align}
     \end{subequations}
     where $C_{\min}^{s_c'}(s') = C_{\min}(s') - C_{\min}(s_c' ) $, $D^{k, s_{d^k}'}_{\min}(s') = D^k_{\min}(s') - D^k_{\min}(s_{d^k}')$, for all $s' \in S$, $ \displaystyle s_c' \in \argmin_{s' \in S} C_{\min}(s') $, $ \displaystyle s_{d^k}' \in \argmin_{s' \in S} D^k_{\min}(s')$, and $\mathfrak{R}_c$, $\mathfrak{R}_{d^k}$,  are defined as 
\begin{align}\label{definition_r_c_and_r_d_k_for_outer_approx}
    &\alpha \mathfrak{R}_c
     = 
    z -  \rho^T \tilde{c}^l   + \alpha \sum_{(s,a) \in \mathcal{K}} \sum_{s' \in S  }  \rho(s,a)  \left(  C_{\min}(s') - C_{\max}(s') \right) \zeta^l  (s' \vert s,a), \nonumber \\
    &\alpha \mathfrak{R}_{d^k}
     =     \xi_k -  \rho^T \tilde{d}^{kl} + \alpha \sum_{(s,a) \in \mathcal{K}} \sum_{s' \in S  }  \rho(s,a)  \left(  D^k_{\min}(s') - D^k_{\max}(s') \right) \zeta^l  (s' \vert s,a), \forall \ k \in \mathbb{K}.
\end{align}
\end{lemma}
\begin{proof}
Let $(z,f)$ be a feasible solution of \eqref{common_JCCMDP_uncertain_TP}. Define 
\[
\rho(s,a) =   \big( (1-\alpha) \gamma^T Q^M_f \allowbreak  \big)(s) f(s,a),  \ \forall \ (s,a) \in \mathcal{K}.
\]
Thus, $\rho $ satisfies \eqref{ICC_random_TPs_outer_approximation_flow_balance_simplified}. As a result,   
\begin{align*}
     z - \tilde{C}_{rc{\text-}tp} \left( \gamma,f \right)
     &
     \leq 
     z -  (1-\alpha) \gamma^T Q^M_f \left(  \tilde{c}_f^l  + 
\alpha  \left( Z_{  f } - Z_{  f }^l \right) C_{\min} + \alpha Z_{  f }^l C_{\max} \right) \\ 
& =
\begin{aligned}[t]
z -  \sum_{(s,a) \in \mathcal{K}} \big( (1-\alpha) \gamma^T Q^M_f \big)(s)  f(s,a)  \bigg\{ \tilde{c}^l(s,a) + \alpha  \sum_{s' \in S }    \big( \zeta  (s' \vert s,a) 
& \\ 
- \zeta^l  (s' \vert s,a) \big) C_{\min}(s') + \alpha \sum_{s' \in S } \zeta^l (s' \vert s,a)  C_{\max}(s')  \bigg\} &
\end{aligned} \\ 
& =
\begin{aligned}[t]
z -  \sum_{(s,a) \in \mathcal{K}}   \rho(s,a)  \bigg\{ \tilde{c}^l(s,a) + \alpha  \sum_{s' \in S }    \big( \zeta  (s' \vert s,a) -  \zeta^l  (s' \vert s,a) \big) C_{\min}(s') 
&
\\
+
\alpha \sum_{s' \in S } \zeta^l (s' \vert s,a)  C_{\max}(s')  \bigg\}
&
\end{aligned} \\ 
& = \alpha \mathfrak{R}_c - \alpha \sum_{(s,a) \in \mathcal{K}} \sum_{s' \in S  } \rho(s,a) C_{\min}(s')  \zeta  (s' \vert s,a)\\ 
& = \alpha \mathfrak{R}_c - \alpha \sum_{(s,a) \in \mathcal{K}} \sum_{s' \in S \backslash s_c' } \rho(s,a) C_{\min}^{s_c'}(s') \zeta  (s' \vert s,a),  
\end{align*}
where the first inequality follows from \eqref{lower_bound_costs_under_random_TPs}, the first equality is the matrix expansion of its preceding expression, the second equality follows from definition of $\rho$, the third equality follows from \eqref{definition_r_c_and_r_d_k_for_outer_approx}, while the last equality follows from the fact that $\displaystyle \sum_{s' \in S} \zeta (s' \vert s, a)(\omega) = 0$,    $(s, a) \in \mathcal{K}$ and $\omega \in \Omega$. Therefore, $(z, \rho)$ is a feasible solution of \eqref{ICC_random_TPs_outer_approximation_objective_simplified}.  
From \eqref{JCC_common_JCCMDP_uncertain_TP}, we have 
\begin{align*}
   \mathbb{P} \left( \tilde{D}^k_{rc{\text-}tp} ( \gamma,f ) \leq \xi_k   \right) \ge \mathbb{P} \left( \tilde{D}^k_{rc{\text-}tp} ( \gamma,f ) \leq \xi_k,  \ \forall  \  k \in \mathbb{K} \right)  \ge p_1,
     \ \forall \ k \in \mathbb{K}.
\end{align*}
Again, using arguments similar as above, we can show that
$(z, \rho)$ is a feasible solution of \eqref{ICC_random_TPs_outer_approximation_constraints_simplified} and the result follows. 
\end{proof}
Since \eqref{ICC_random_TPs_outer_approximation_simplified} is an ILCCP problem, 
we apply  \Cref{general lower bound} to obtain a lower bound to the optimal value of \eqref{common_JCCMDP_uncertain_TP}.  
\begin{theorem}
Let conditions \ref{non_negativity_perturb_TP_assumption} and \ref{costs_known_TP_assumption} of \Cref{existence_of_perturbations} hold. A lower bound  to the optimal value of \eqref{common_JCCMDP_uncertain_TP}   
 is given by the optimal value of the following LP problem
\begin{align} 
    \min_{z, \rho, \mathfrak{m}^{rp}_c, \big( \mathfrak{m}^{rp}_{d^k} \big)_{k \in \mathbb{K}}} & \ z \nonumber \\
 \textnormal{s.t. } 
 -\mathfrak{R}_c \leq (1 - & p_0) \mathfrak{m}^{rp}_c,
 \ \lambda^{rp}_c \leq \mathfrak{m}^{rp}_c, \ -\mathfrak{R}_{d^k} \leq (1 - p_1) \mathfrak{m}^{rp}_{d^k},
 \ \lambda^{rp}_{d^k} \leq \mathfrak{m}^{rp}_{d^k}, \ \forall \ k \in \mathbb{K}, \nonumber\\
  \sum_{(s,a) \in \mathcal{K}}  \sum_{s' \in S \backslash s_c'  } &\ \rho(s,a) C_{\min}^{s_c'}(s') \zeta^u (s' \vert s,a)
- \mathfrak{R}_c  \leq \mathfrak{m}^{rp}_c,
\nonumber
\\
 \sum_{(s,a) \in \mathcal{K}} \sum_{s' \in S \backslash s_c'  } & \  \rho(s,a) C_{\min}^{s_c'}(s') \zeta^l (s' \vert s,a)  \leq \mathfrak{R}_c, \nonumber \\
 \ \sum_{(s,a) \in \mathcal{K}} \sum_{s' \in S \backslash s_{d^k}'  } & \   \rho(s,a) D^{k, s_{d^k}'}_{\min}(s') \zeta^u (s' \vert s,a)     
- \mathfrak{R}_{d^k} \leq \mathfrak{m}^{rp}_{d^k}, \ \forall \ k \in \mathbb{K}, \nonumber
\\
 \sum_{(s,a) \in \mathcal{K}} \sum_{s' \in S \backslash s_{d^k}'  } & \   \rho(s,a) D^{k, s_{d^k}'}_{\min}(s') \zeta^l (s' \vert s,a)  \leq \mathfrak{R}_{d^k},  \, \forall \, k \in \mathbb{K}, \nonumber  \\
 & \hspace{-1.5cm} \rho \in \mathcal{Q}_M^{\alpha}(\gamma), \label{ICC_outer_approximation_using_rho}
\end{align}
where  $\lambda^{rp}_c >0$ , $\lambda^{rp}_{d^k} >0$,   are constants, and  $\mathfrak{R}_c$, $\mathfrak{R}_{d^k}$,  $k \in \mathbb{K}$, are defined as in \eqref{definition_r_c_and_r_d_k_for_outer_approx}. 
\end{theorem}
\begin{proof}
 From \Cref{lower-bd-lemma}, it is enough to show that for every feasible solution of \eqref{ICC_random_TPs_outer_approximation_simplified}, there exists a feasible solution of \eqref{ICC_outer_approximation_using_rho}. 
  Let $(z, \rho)$ be a feasible solution of \eqref{ICC_random_TPs_outer_approximation_simplified}. Define 
\begin{align*}
\mathfrak{m}^{rp}_c & = \max\bigg( \lambda^{rp}_c, \ \sum_{(s,a) \in \mathcal{K}} \sum_{s' \in S \backslash s_c' }  \rho(s,a) C_{\min}^{s_c'}(s') \zeta^u (s' \vert s,a)  - \mathfrak{R}_c  \bigg),  \\ 
\mathfrak{m}^{rp}_{d^k} & = 
\begin{aligned}[t]
\max\bigg( \lambda^{rp}_{d^k}, \ \sum_{(s,a) \in \mathcal{K}} \sum_{s' \in S \backslash s_{d^k}' }  \rho(s,a) D^{k, s_{d^k}'}_{\min}(s')  \zeta^u (s' \vert s,a)  - \mathfrak{R}_{d^k}  \bigg),  \ \forall \ k \in \mathbb{K}.
\end{aligned}
\end{align*}
Since the coefficient of each random variable in the chance constraints of \eqref{ICC_random_TPs_outer_approximation_simplified} is non-negative, the result follows directly from \Cref{general lower bound}.
\end{proof}
\VV{
\begin{remark}
    In terms of state-action pairs, the time required to solve the above LP problem can be bounded above by $O( \vert \mathcal{K} \vert^{\frac{1}{2}} \vert \ln \hat{\epsilon} \vert )$ \cite{wright1997primal}. Assuming same action set $A$ at each state, we obtain $O( \vert S \vert^{\frac{1}{2}} \vert A \vert^{\frac{1}{2}} \vert  \ln \hat{\epsilon} \vert )$ \VS{as the upper bound on time}. 
\end{remark}
}
\subsection{Bounds on gap}\label{Theoretical guarantees TPs}
Similar to Section \ref{Theoretical guarantees costs}, 
we derive  extremal bounds \text{UB}$_{rc{\text-}tp}^{(u)}$ and \text{LB}$_{rc{\text-}tp}^{(l)}$ on the optimal values of the approximations from  Sections \ref{Inner approximations random TPs} and \ref{Outer approximations random TPs}, respectively, and summarize them in 
 \Cref{Theoretical_Guarantees_TPs}. The Gap(\%), defined by \eqref{gap_per}, in this case belongs to interval $[0, G_{rc\text{-}tp}]$, where $G_{rc\text{-}tp}= \displaystyle \frac{ \text{UB}_{rc\text{-}tp}^{(u)} \text{ - } \text{LB}_{rc\text{-}tp}^{(l)} }{ \text{LB}_{rc\text{-}tp}^{(l)}}\times  100$. In Table \ref{Theoretical_Guarantees_TPs}, 
$\mathfrak{g}(\rho) = \displaystyle \alpha \sum_{(s,a) \in \mathcal{K}} \sum_{s' \in S  } \rho(s,a)  \left(  C_{\max}(s') - C_{\min}(s') \right) \zeta^l  (s' \vert s,a) $. In the numerical experiments, we observe a significant reduction in the gap  obtained
by solving our approximations as compared to the theoretical ones, as in the case of random costs.
\begingroup
    \small
\renewcommand\arraystretch{1.6}
    \begin{longtable}{c p{13cm}}
    \caption{Extremal bounds.}\\
    \toprule
    \endfirsthead
\caption{Extremal bounds.}\\
\toprule
\endhead
\midrule
\multicolumn{2}{r}{\footnotesize\itshape Continued on the next page}
\endfoot
\endlastfoot
    Upper bounds & \multicolumn{1}{c}{\text{UB}$_{rc{\text-}tp}^{(u)}$} \\ 
    \midrule
    \eqref{C_H_Sub_G_random_TPs_convex_approximation} & \\ 
    \cline{1-1}
 Part \ref{One-sided Chebyshev inequality random TPs convex approximation} 
    & $\begin{aligned}[t]
    \max_{\rho \in \mathcal{Q}^{\alpha}(\gamma)}  \rho^T \tilde{c}^u  -        \mathfrak{g}(\rho) 
    + \alpha \sqrt{\frac{p_0}{1 - p_0}} \, \bigg(\sum_{ (s,a) \in \mathcal{K} }  \rho (s,a)  \big\| \big(\Sigma_{C \zeta} ^{\frac{1}{2}}  \big)_{(s,a)} \big\|_2 \bigg) 
    \end{aligned}$
    \\
    Part \ref{Hoeffding inequality random TPs convex approximation} 
    & 
    $\begin{aligned}[t]
    \max_{\rho \in \mathcal{Q}^{\alpha}(\gamma)}  \rho^T \tilde{c}^u  -  \mathfrak{g}(\rho)
    + \alpha \sqrt{-\frac{1}{2} \ln   (1-p_0)  }  \bigg( \sum_{ (s,a) \in \mathcal{K} }  
 \displaystyle   \sum_{s' \in S \backslash s_c} \rho (s,a) C_{\max}^{s_c}(s')  (  \zeta^u(s' \vert s, \\ a)  
  \vphantom{\sum_{s' \in S \backslash s_c}  }
 - \zeta^l (s' \vert s,a) ) \bigg) 
    \end{aligned}$
    \\
   \cline{1-1}
   \eqref{inner_approximation_using_bernstein_inequality_theorem_TP} &
   $\begin{aligned} [t]
 \max_{\rho \in \mathcal{Q}^{\alpha}(\gamma)}  \rho^T \tilde{c}^u - \mathfrak{g}(\rho)
+   
 \alpha \sum_{(s,a) \in \mathcal{K}}  \sum_{s' \in S \backslash s_c} \rho(s,a) C_{\max}^{s_c}(s')      \zeta^u  (s' \vert s,a)         
  - \frac{\alpha}{h_0^{rp}} \ln(1-p_0) &
  \end{aligned}$
   \\
   \midrule
   Lower bound & \multicolumn{1}{c}{\text{LB}$_{rc{\text-}tp}^{(l)}$} \\ 
   \midrule
    \eqref{ICC_outer_approximation_using_rho} &  
    $\begin{aligned}[t]
        \min_{\rho \in \mathcal{Q}^{\alpha}(\gamma)}  
\rho^T \tilde{c}^l 
+ 
  \mathfrak{g}(\rho)
+ 
  \alpha  \sum_{(s,a) \in \mathcal{K}} \sum_{s' \in S \backslash s_c'  }  \  \rho(s,a) C_{\min}^{s_c'}(s') \zeta^l (s' \vert s,a)
    \end{aligned}$
   \\
   \bottomrule  \label{Theoretical_Guarantees_TPs}
    \end{longtable}
    \endgroup
    \section{Numerical experiments} \label{Num}
    We perform all the numerical experiments using  CVX package in MATLAB optimization toolbox,  on an Intel(R) 64-bit Core(TM) i5-1240P CPU @ 1.70GHz with 16.0 GB RAM machine. For the JCCMDP problem under random running costs, we perform numerical experiments on randomly generated instances of a queueing control problem studied in \cite{varagapriya2022joint}.  For the JCCMDP problem under both random running costs and  transition probabilities, we consider randomly generated CMDPs belonging to a well-known class of Markov decision problems called Garnets \cite{Archibald_1995, SBIRL, varagapriya2022joint}. In both cases, we solve the convex approximations and analyse the gaps between upper and lower bounds obtained by varying parameters and compare them with the extremal bounds derived in Tables \ref{Theoretical_Guarantees_costs} and \ref{Theoretical_Guarantees_TPs}. In all our numerical experiments, we fix $p_0 = 0.9$ and  $p_1 = 0.9$.
\subsection{Queueing control problem under random running costs}\label{Queueing control problem}
We consider a single queue system with service and admission controllers \cite{altman, varagapriya2022joint}. The maximum queue length $L$ is considered to be finite and it represents the number of states. The state space is given by $S=\{0,1,2,\ldots, L\}$, where $s=0$ denotes an empty queue, while $s = L$ denotes a full queue. A service controller finishes the service with a given probability level $a^1$ and an admission controller admits a new customer with a given probability level $a^2$. It is assumed that $a^1$ and $a^2$   can take finitely many values from the sets $A^1\subseteq [a_{\min}^1, a_{\max}^1]$ and $A^2\subseteq [a_{\min}^2,a_{\max}^2]$, respectively,  where $0<a_{\min}^1\leq a_{\max}^1<1$ and $0\leq a_{\min}^2\leq a_{\max}^2<1$. Since the service and admission controllers have a common objective, they are considered as a single decision-maker, and their joint action is represented as a tuple $(a^1,a^2)\in A^1\times A^2$.  
  It is assumed that  admission does not take place when $s = L$.
   The transition probabilities are given as in \cite{altman}
   \begin{align*}
p\big( s'\vert s,(a^1,a^2) \big)=\begin{cases}
a^1(1-a^2) & 1\leq s\leq L-1, s'=s-1,\\
a^1a^2+(1-a^1)(1-a^2) & 1\leq s\leq L-1, s'=s,\\
(1-a^1)a^2 & 0\leq s\leq L-1, s'=s+1,\\
1-(1-a^1)a^2 & s'=s=0,\\
1-a^1 & s=L,s'=L,\\
a^1 & s=L,s'=L-1.
\end{cases}
\end{align*}
For every  state-action pair $(s,(a^1,a^2))$, three different costs, namely, holding cost $\tilde{c}(s)$, service cost $\tilde{d}^1(a^1)$, and no-admission cost $\tilde{d}^2(a^2)$ are incurred. The aim is to minimize the expected discounted holding cost subject to upper bounds on the expected discounted service and no-admission costs. Since these costs are not exactly known in advance, they are better modeled using random variables without imposing a requirement that their exact distribution is known. We examine how the gaps between the upper and lower bounds vary with different approximations. 
We use the convex approximations derived in 
parts \ref{One-sided Chebyshev inequality random costs convex approximation}-\ref{Hoeffding inequality random costs convex approximation} of \Cref{C_H_Sub_G_random_costs_convex_approximation_statement}  and \Cref{Bernstein inequality random costs} to obtain upper bounds and \Cref{linear_approximation_random_costs_statement} to obtain a lower bound to the optimal value of the   JCCMDP problem \eqref{ICCP_random_costs}. These bounds are obtained assuming that the mean, variance, and upper and lower bounds of the random running costs are known. However, to assess how these gaps perform when the exact distribution of the random costs is known, we  solve the approximations proposed in Section 3.2 of \cite{varagapriya2022joint} with various distributions and compare all the resulting gaps. In addition, we focus on the upper bound approximations and analyse how they vary as the availability of information on the distribution of the running costs improves.

\begin{figure*}[hbt!]
    \centering
    \begin{subfigure}[c]{0.5\linewidth}
        \centering
        \includegraphics[scale=0.63]{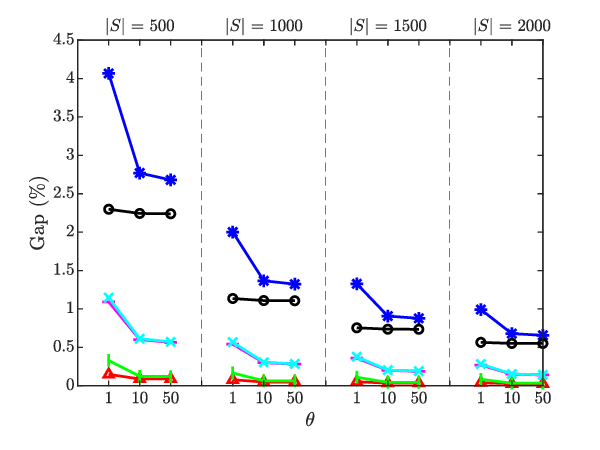}
        \vspace*{-3mm}
    \caption{\footnotesize  }
    \label{random_running_costs_optimal_values_plot}
    \end{subfigure}%
    \begin{subfigure}[c]{0.5\textwidth} 
       \centering
        \includegraphics[scale=0.63]{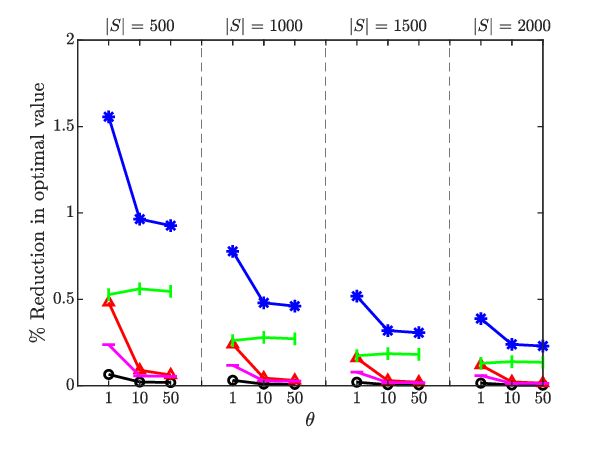}
         \vspace*{-3mm}
    \caption{\footnotesize   }
    \label{Percent_reduction_UB_values_plot}
    \end{subfigure}
   \vspace*{-5mm}
   \caption{ \footnotesize In \ref{random_running_costs_optimal_values_plot}, $`\ast$',\textnormal{`o'}, and  $`\Delta$',  $`\vert$', $`-$', and $`\times$' indicate the average gaps  corresponding to \eqref{C_H_Sub_G_random_costs_convex_approximation} under part \ref{One-sided Chebyshev inequality random costs convex approximation} of \Cref{C_H_Sub_G_random_costs_convex_approximation_statement}, \eqref{Bernstein_inequality_random_costs},  and 
   approximations of \cite{varagapriya2022joint} under Gaussian, T-Student(5),  Logistic, and T-Student(1) distributions, respectively. In \ref{Percent_reduction_UB_values_plot}, $`\ast$',\textnormal{`o'},  $`\Delta$', $`\vert$', and $`-$' indicate the percentage reduction in the average upper bound values when  we switch the approximations in the following  order:
\eqref{C_H_Sub_G_random_costs_convex_approximation} under part \ref{One-sided Chebyshev inequality random costs convex approximation} of  \Cref{C_H_Sub_G_random_costs_convex_approximation_statement} $\rightarrow$  approximations of \cite{varagapriya2022joint} under  T-Student(1) $\rightarrow$ Logistic $\rightarrow$  \eqref{C_H_Sub_G_random_costs_convex_approximation} under 
part \ref{sub_G_theorem_statement_random_costs} of \Cref{C_H_Sub_G_random_costs_convex_approximation_statement} $\rightarrow$ approximations of \cite{varagapriya2022joint} under T-Student(5) $\rightarrow$ Gaussian distributions, respectively.
   }
\end{figure*}

As similar to \cite{varagapriya2022joint}, we fix $\alpha = 0.9$,  $ \allowbreak A^1 \times A^2=\{(0.2, 0.75, 0.9),(0, \allowbreak 0.5,  0.8)\}$, and  uniformly generate $\gamma$.  The components of the mean vectors $\mu_{\tilde{c}}$, $\mu_{\tilde{d}^1}$, and $\mu_{\tilde{d}^2}$ are defined as  $\mu_{\tilde{c}}(s)=s$, $\mu_{\tilde{d}^1}(a^1)=3(1+a^1)^2$, and $\mu_{\tilde{d}^2}(a^2)=10-3a^2$,  $s \in S$, $(a^1, a^2) \in A^1 \times A^2$. We take $\xi_1= 11.30$ and $\xi_2=11.35$.  For each running cost vector, we randomly generate diagonal covariance matrices, $\Sigma_{\tilde{c}}$, $\Sigma_{\tilde{d}^1}$, and $\Sigma_{\tilde{d}^2}$ with values in $(0,0.8)$ and for each $(s, a) \in \mathcal{K}$, we fix their $(s,a)^{th}$ diagonal entries  as $\sigma^2_{\tilde{c}} (s,a)$, $\sigma^2_{\tilde{d}^1} (s,a)$, and $\sigma^2_{\tilde{d}^2} (s,a)$, respectively. Following the procedure specified in \cite{peng2022bounds}, we generate 3000 samples from Gaussian distribution with mean vector $\mu_{\tilde{c}}$ and covariance matrix $\Sigma_{\tilde{c}}$. We fix the maximum and minimum values of these samples as $\tilde{c}^u$ and $\tilde{c}^l$, respectively. Similarly, the upper and lower bounds of other random running costs are fixed.  
 In \eqref{Bernstein_inequality_random_costs}, we fix $h_0^{rc} = 10$, $h_1^{rc} = 10$, and $h_2^{rc} = 10$.  In \eqref{linear_approximation_random_costs}, we fix $\lambda^{rc}_c = 10^{-5} $, $\lambda^{rc}_d = 10^{-5}$, $N = 20$, and generate the points $ y_1^i,  y_2^i $ in $(0.1,1)$, for $i = 1, 2, \ldots, 20$. For the approximations of \cite{varagapriya2022joint}, we assume four elliptical distributions for random running costs, namely, Gaussian,   T-Student with degrees of freedom 1 and 5  (denoted henceforth by T-Student(1) and T-Student(5), respectively), and Logistic distributions \cite{fang}, with the previously specified mean vectors and covariance matrices as parameters. We  fix $N=20$ to approximate the quantile function (see Proposition 2 of \cite{varagapriya2022joint}). We vary the values of $\theta $ and $\vert S \vert$  and generate 50 random instances for each problem. We compute the average of the gaps between upper and lower bounds constructed in \Cref{JCCMDP problem under random running costs} and those constructed in \cite{varagapriya2022joint}, where  the gap between the upper and lower bounds is defined as in \eqref{gap_per}.  
We summarize the average gaps obtained for different approximations in \Cref{random_running_costs_optimal_values_plot}. 
 \noindent
We observe a similar pattern across $\vert S \vert$; the average gaps corresponding to all the approximations reduce with $\theta$ and $\vert S \vert$ with the best gap obtained corresponding to the approximation of \cite{varagapriya2022joint} under Gaussian distribution. While not practical, this indicates that knowing the exact distribution is beneficial.
Among the approximations \eqref{C_H_Sub_G_random_costs_convex_approximation} under 
parts \ref{One-sided Chebyshev inequality random costs convex approximation}-\ref{Hoeffding inequality} of \Cref{C_H_Sub_G_random_costs_convex_approximation_statement} and  \eqref{Bernstein_inequality_random_costs}, 
  we obtain the best gap corresponding to \eqref{Bernstein_inequality_random_costs}. We highlight that out of the 50 instances we ran, at most 3 instances were feasible corresponding to \eqref{C_H_Sub_G_random_costs_convex_approximation} under part \ref{Hoeffding inequality random costs convex approximation} of \Cref{C_H_Sub_G_random_costs_convex_approximation_statement}, for each $\theta$ and $\vert S \vert$, and hence, we do not plot the gap corresponding to it.   Since every instance of \eqref{C_H_Sub_G_random_costs_convex_approximation} under part \ref{One-sided Chebyshev inequality random costs convex approximation} of \Cref{C_H_Sub_G_random_costs_convex_approximation_statement} is feasible, it follows from \Cref{comparison_f_k_values} that the infeasibility of most of the instances under part \ref{Hoeffding inequality random costs convex approximation} is due to a high value of $V_k^*$, which depends on the upper and lower bounds of running costs. Moreover, the  approximation of part \ref{sub_G_theorem_statement_random_costs} is better than the approximation of part \ref{One-sided Chebyshev inequality random costs convex approximation}.
 This can be observed from  \Cref{Percent_reduction_UB_values_plot}.   
 With the previously specified mean vectors and covariance matrices as parameters, the average upper bound values of \eqref{ICCP_random_costs} reduce as we move from the optimal values of \eqref{C_H_Sub_G_random_costs_convex_approximation} under 
part \ref{One-sided Chebyshev inequality random costs convex approximation} to part \ref{sub_G_theorem_statement_random_costs} of \Cref{C_H_Sub_G_random_costs_convex_approximation_statement} and subsequently to the approximation of \cite{varagapriya2022joint} under a Gaussian distribution, i.e.,    as we gain more insight into the distribution of the random running costs. In addition, we observe reductions in upper bound values as the tails of the elliptical distributions become lighter and for a given $\theta $, the percentage of reduction decreases with $\vert S \vert$ at each switch.
  In terms of the average CPU time taken, among the upper bound approximations, it takes maximum time to solve \eqref{Bernstein_inequality_random_costs},  reaching up to 46.8388
seconds when $(\theta, \vert S \vert) = (50, 2000)$, while the minimum time varies among approximations, with the lowest being 1.0472 seconds corresponding to \eqref{C_H_Sub_G_random_costs_convex_approximation} under part \ref{One-sided Chebyshev inequality random costs convex approximation} of \Cref{C_H_Sub_G_random_costs_convex_approximation_statement} when $(\theta, \vert S \vert) = (10, 500)$. On the other hand, among the lower bound approximations, it takes lesser time to solve \eqref{linear_approximation_random_costs} than the approximations in \cite{varagapriya2022joint}, with the time ranging between 0.7288  and 15.1175 
seconds. 

Furthermore, in \Cref{Theoretical_guarantees_costs_graph}, we summarize the average of the percentage of reduction in the gaps when transitioning from the extremal bounds given in \Cref{Theoretical_Guarantees_costs} to the approximations from Sections \ref{Inner approximations for random costs} and \ref{Outer approximation for random costs} along with similar derivations for the approximations of \cite{varagapriya2022joint}.
We observe that the percentage of reduction increases with   $\theta$, for each $\vert S \vert$.

\subsection{Randomly generated CMDPs under random running costs and  transition probabilities}\label{CMDP_with_random_TPs}
We generate random finite CMDPs as studied in Section 4.2 of \cite{varagapriya2022joint} and assume that the running costs and  transition probabilities in the system are random. We construct the system as described in \cite{SBIRL, varagapriya2022joint}, by defining a tuple $(\vert S \vert, \vert A \vert, \vert B_F \vert)$, where the branching factor, $ \vert B_F \vert$, denotes the number of states reachable from every state-action pair.  Suppose, from $(s,a) \in \mathcal{K}$, $(s^i)_{i=1}^{\vert B_F \vert}$ states are reachable. We fix the associated expected transition probabilities by randomly generating $\vert B_F \vert-1$ values from $(0,1)$, denoted by, $(q^i)_{i=1}^{\vert B_F \vert-1}$, using the function \textbf{sort(rand($\vert B_F \vert-1,1 $))}. We define the expected transition probabilities by
\begin{align*}
\mu(s'\vert s,a)=\begin{cases}
q^i-q^{i-1} &   s'=s^i,\, i=1,2,\ldots,\vert B_F \vert,\\
0 & \text{otherwise},
\end{cases}
\end{align*}
where $q^0 = 0$ and $q^{\vert B_F \vert } = 1$.
 We fix  $\vert S \vert = 250$, $\vert A \vert = 15$, $\vert B_F \vert = 125$,  $K =15$, and  uniformly generate $\gamma$.  We fix the values of $\tilde{c}^u$ and $\tilde{c}^l$ as discussed in Section \ref{Queueing control problem} using mean vector and diagonal covariance matrix randomly generated from the intervals $(50,70)$ and $(0,0.4)$, respectively. The values of $\tilde{d}^{ku}$ and $\tilde{d}^{kl}$, for all $k = 1,2, \ldots,15$, are similarly fixed with mean vectors and diagonal covariance matrices randomly generated from the intervals $(50,100)$ and $(0,0.4)$, respectively.  We randomly generate   $\xi_k \in \mathbb{R}^{15}$, from the interval    $(80, 90)$.   For each $(s, a) \in \mathcal{K} $ and $s' \in S$,  we 
 generate  3000 samples from uniform  distribution from the interval $\big( -\eta \mu(s' \vert s, a),\eta( 1- \mu(s' \vert s, a) ) \big)$. We fix the maximum and minimum values
of these samples as $\zeta^u (s' \vert s, a)$ and $\zeta^l (s' \vert s,a)$, respectively. Whenever $\mu(s' \vert s, a) = 0$, we fix $\zeta^l (s' \vert s, a)=0$.   
   In \eqref{inner_approximation_using_bernstein_inequality_theorem_TP}, we fix $h_0^{rp} = 10,$ $h_k^{rp} = 10$, while in \eqref{ICC_outer_approximation_using_rho}, we fix $\lambda_c^{rp} = 10^{-5},$ $\lambda_{d^k}^{rp} = 10^{-5}$,   $k = 1, 2, \ldots, 15$.  We vary the values of $\eta$ and $\alpha$ and generate 50 random instances for each problem. We compute the average of the gaps between upper and
 lower bounds constructed in \Cref{JCCMDP problem under random running costs and transition probabilities} 
 assuming that the mean, variance, and upper and lower bounds of the random perturbations are known.

\begin{figure*}[hbt!]
    \centering
    \begin{subfigure}[c]{0.5\linewidth}
        \centering
    \includegraphics[scale=0.63]{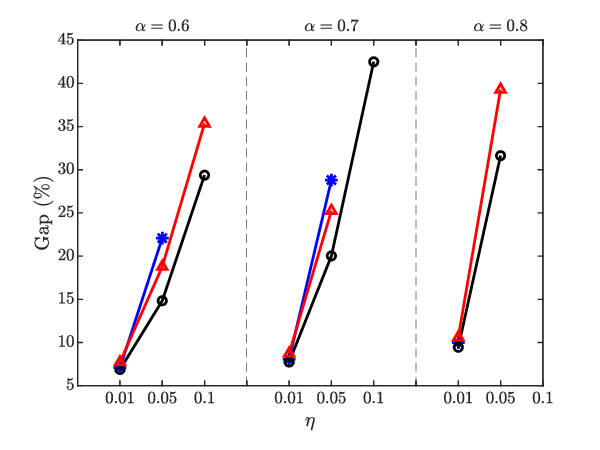}
     \vspace*{-3mm}
    \caption{\footnotesize  }
\label{random_transition_probabilities_optimal_values_plot}
 \end{subfigure}%
    \begin{subfigure}[c]{0.5\textwidth} 
    \centering
      \includegraphics[scale=0.63]{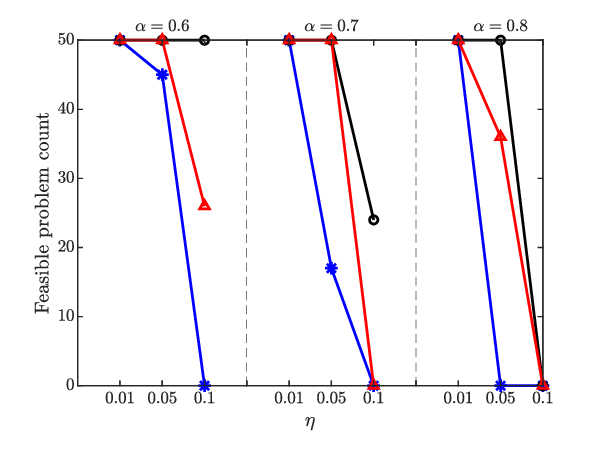}
         \vspace*{-3mm}
    \caption{\footnotesize   }
    \label{Feasible_solution_count_TPs_plot}
    \end{subfigure}
   \vspace*{-5mm}
    \caption{ \footnotesize In \ref{random_transition_probabilities_optimal_values_plot}, $`\ast$', $`\Delta$', and \textnormal{`o'}  indicate the average gaps  corresponding to \eqref{C_H_Sub_G_random_TPs_convex_approximation} under parts \ref{One-sided Chebyshev inequality random TPs convex approximation}-\ref{Hoeffding inequality random TPs convex approximation} of \Cref{C_H_Sub_G_random_TPs_convex_approximation_thm_statement}, and \eqref{inner_approximation_using_bernstein_inequality_theorem_TP}, respectively. In \ref{Feasible_solution_count_TPs_plot}, these indicate the number of instances with a feasible solution. 
   }
\end{figure*}

\begin{figure*}[hbt!]
    \centering
    \begin{subfigure}[c]{0.5\linewidth}
        \centering
    \includegraphics[scale=0.63]{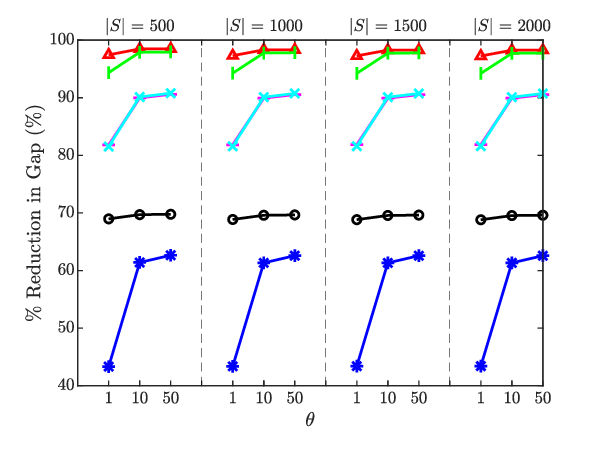}
     \vspace*{-3mm}
    \caption{\footnotesize  }
\label{Theoretical_guarantees_costs_graph}
 \end{subfigure}%
    \begin{subfigure}[c]{0.5\textwidth} 
    \centering
      \includegraphics[scale=0.63]{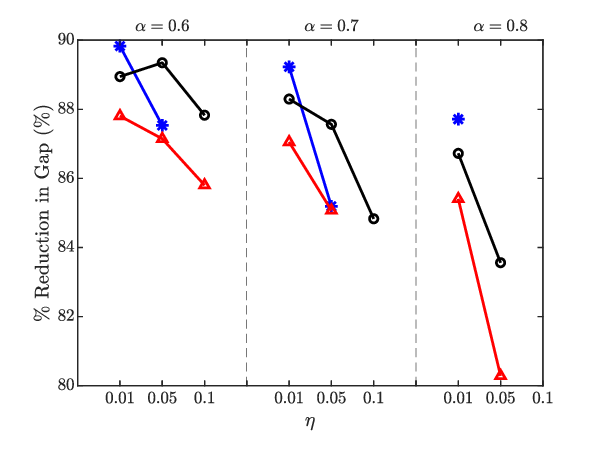}
         \vspace*{-3mm}
    \caption{\footnotesize   }
    \label{Theoretical_guarantees_TPs_graph}
    \end{subfigure}
   \vspace*{-5mm}
    \caption{ \footnotesize In \ref{Theoretical_guarantees_costs_graph}, $`\ast$',\textnormal{`o'}, and $`\Delta$',  $`\vert$', $`-$', and $`\times$'  indicate the average percentage reduction in the gaps corresponding to \eqref{C_H_Sub_G_random_costs_convex_approximation} under part \ref{One-sided Chebyshev inequality random costs convex approximation} of \Cref{C_H_Sub_G_random_costs_convex_approximation_statement}, \eqref{Bernstein_inequality_random_costs},  and approximations of \cite{varagapriya2022joint} under Gaussian, T-Student(5),  Logistic, and T-Student(1) distributions, respectively. In \ref{Theoretical_guarantees_TPs_graph}, $`\ast$',\textnormal{`o'}, and $`\Delta$'  correspond to \eqref{C_H_Sub_G_random_TPs_convex_approximation} under part \ref{One-sided Chebyshev inequality random TPs convex approximation}  of \Cref{C_H_Sub_G_random_TPs_convex_approximation_thm_statement}, \eqref{inner_approximation_using_bernstein_inequality_theorem_TP}, and  \eqref{C_H_Sub_G_random_TPs_convex_approximation} under part \ref{Hoeffding inequality random TPs convex approximation}  of \Cref{C_H_Sub_G_random_TPs_convex_approximation_thm_statement}, respectively.
   }
\end{figure*}

 We summarize the number of instances giving feasible solutions in \Cref{Feasible_solution_count_TPs_plot}. 
 For upper bound approximations providing feasible solutions in at least 15 instances, the corresponding average gaps are summarized in \Cref{random_transition_probabilities_optimal_values_plot}.
  We observe that the gaps increase with $\eta$ and $\alpha$, with the best gap obtained corresponding to \eqref{inner_approximation_using_bernstein_inequality_theorem_TP}.   For relatively higher values of $\eta$ and $\alpha$,  fewer than 50 instances of the upper bound approximations yield feasible solutions, and the resulting gaps are larger than those that provide feasible solutions in every instance.    
In terms
of the average CPU time taken, among the upper bound approximations,  it takes  maximum time
 to solve \eqref{inner_approximation_using_bernstein_inequality_theorem_TP}  reaching up to 64.8853 seconds when $(\eta, \alpha) = (0.01, 0.8)$ and minimum time
 to solve \eqref{C_H_Sub_G_random_TPs_convex_approximation} under part \ref{Hoeffding inequality random TPs convex approximation} of \Cref{C_H_Sub_G_random_TPs_convex_approximation_thm_statement} with the lowest being 4.6925
 seconds when $(\eta, \alpha) = (0.01, 0.6)$. On the other hand, the time to solve \eqref{ICC_outer_approximation_using_rho} ranges between 0.9294 and 2.4866 seconds.  
Furthermore, in \Cref{Theoretical_guarantees_TPs_graph}, we summarize the average of the percentage of reduction in the gaps when transitioning from the extremal bounds given in \Cref{Theoretical_Guarantees_TPs} to the approximations from Sections \ref{Inner approximations random TPs} and \ref{Outer approximations random TPs}. We observe that the percentage of reduction decreases with  $\alpha$, for each $\eta$.

\VV{Our numerical experiments confirm that the convex approximations proposed in \Crefrange{JCCMDP problem under random running costs}{JCCMDP problem under random running costs and transition probabilities} scale well and remain computationally efficient across problem sizes.  Notably, the SOCP and LP problems are efficiently solvable due to their underlying polynomial-time complexity.} \AL{Notice that if the number of states and/or the number of actions increase significantly, the computational efforts will increase accordingly.}
\section{Conclusions}
\label{sec:conclusions}
We study a JCCMDP problem where either running costs are random or both running costs and transition probabilities are random.
When the uncertainty is present only in the running cost vectors, we assume that the dependency among the random constraint vectors is driven by a Gumbel-Hougaard copula, whereas for the case when both running costs and transition probabilities are random, we make no assumption on the dependency structure of underlying random parameters. In both cases, we use standard probability inequalities, including one-sided Chebyshev, Bernstein, and Hoeffding inequalities, to derive upper bound approximations. Further, we show that a lower bound can be obtained from the optimal solution of an LP  problem. 
To analyse the gap between upper and lower bounds, we perform numerical experiments on randomly generated instances of a queueing control problem and Garnets. We compare our results with the existing approximations for the case of elliptical distributions  \cite{varagapriya2022joint}. 
We observe from the numerical experiments that an increased knowledge about the 
 distribution is consistently advantageous and leads to better upper
 bounds.
In addition, among the approximations proposed in this paper, those derived using Bernstein's inequality provide the best gap in most of the instances. 
For the case when both running costs and transition probabilities are random, we restrict the problem to the class of stationary policies. \VV{Future research could explore sufficient conditions that guarantee the optimality of stationary policies, as well as derive tractable convex approximations for the JCCMDP problem when the true distribution of underlying parameters is known to belong to a predetermined uncertainty set. Notably, our current approximations assume exact knowledge of parameters such as the mean, variance, and bounds of the random vectors--an assumption that may not hold in practice. In such cases, a promising direction would be to leverage \VS{ robust optimization formulations which consider the worst-case approach with respect to underlying parameters belonging to a predetermined uncertainty set.}}
\section*{Acknowledgments}
The research of first author was supported by CSIR, India.
The research of second  author was supported by DST/CEFIPRA Project No.
IFC/4117/DST-CNRS-5th call/2017-18/2. 
The research of third author was supported by the French government under the France 2030 program, reference ANR-11-IDEX-0003 within the OI H-Code.

 \printbibliography
\end{document}